\documentclass[a4paper,11pt,reqno,noindent]{amsart}
\usepackage[centertags]{amsmath}
\usepackage{amsfonts,amssymb,amsthm,dsfont,cases,amscd,esint,enumerate}
\usepackage[T1]{fontenc}
\usepackage[english]{babel}
\usepackage[applemac]{inputenc}
\usepackage{newlfont}
\usepackage{color}
\usepackage[body={15cm,21.5cm},centering]{geometry} 
\usepackage{fancyhdr}
\usepackage{enumerate}

\theoremstyle{plain}
\newtheorem{theor10}{Theorem}

\newtheorem{prop10}{Proposition}

\newtheorem{cor10}{Corollary}

\newtheorem{theor0}{Theorem}[section]
\newenvironment{theor}
  {\pushQED{\qed}\begin{theor0}}
  {\popQED\end{theor0}}
\newtheorem{lem0}[theor0]{Lemma}
\newenvironment{lem}
  {\pushQED{\qed}\begin{lem0}}
  {\popQED\end{lem0}}
\newtheorem{prop0}[theor0]{Proposition}
\newenvironment{prop}
  {\pushQED{\qed}\begin{prop0}}
  {\popQED\end{prop0}}
\newtheorem{cor0}[theor0]{Corollary}
\newenvironment{cor}
  {\pushQED{\qed}\begin{cor0}}
  {\popQED\end{cor0}}
\newtheorem{propr0}[theor0]{Property}

\newtheorem{hyp0}[theor0]{Hypothesis}

\newtheorem{result0}[theor0]{Result}

\newtheorem{conj0}[theor0]{Conjecture}

\newtheorem{heur0}[theor0]{Heuristics}

\theoremstyle{definition}
\newtheorem{defin0}[theor0]{Definition}
\newenvironment{defin}
  {\pushQED{\qed}\begin{defin0}}
  {\popQED\end{defin0}}
\newtheorem{rems0}[theor0]{Remarks}

\newtheorem{ex0}[theor0]{Example}

\newtheorem{exs0}[theor0]{Examples}

\newtheorem{rem0}[theor0]{Remark}
\newenvironment{rem}
  {\pushQED{\qed}\begin{rem0}}
  {\popQED\end{rem0}}
\newtheorem{qu0}[theor0]{Question}

\newtheorem{qus0}[theor0]{Questions}

  \newtheorem{as0}[theor0]{Assumption}

\mathchardef\emptyset="001F
\numberwithin{equation}{section}

\newcommand{\nablas}{\nabla_{\operatorname{sym}}}

\newcommand{\N}{\mathbb N}

\newcommand{\elec}{\mathrm{elec}}

\newcommand{\calB}{\mathcal{B}}
\newcommand{\calP}{\mathcal{P}}
\newcommand{\calM}{\mathcal{M}}

\newcommand{\Hc}{\mathcal{H}}

\newcommand{\calS}{\mathcal S}

\newcommand{\R}{\mathbb R}
\newcommand{\Z}{\mathbb Z}

\newcommand{\m}{\mathrm{min}}

\newcommand{\loc}{{\operatorname{loc}}}
\newcommand{\Id}{\operatorname{Id}}

\newcommand{\supp}{\operatorname{supp}}

\newcommand{\dist}{\operatorname{dist}}

\newcommand{\uloc}{{\operatorname{uloc}}}

\newcommand{\step}[1]{\noindent \textit{Step} #1.}
\newcommand{\substep}[1]{\noindent \textit{Substep} #1.}

\newcommand{\pr}[1]{\mathbb{P}\left[ #1 \right]}

\newcommand{\expec}[1]{\mathbb{E}\left[ #1 \right]}

\newcommand{\var}[1]{\mathrm{Var}\left[#1\right]}

\newcommand{\cov}[2]{\operatorname{Cov}\left[{#1};{#2}\right]}

\newcommand{\oscd}[2]{\operatorname{osc}_{#2}^2{#1}}
\newcommand{\osc}[2]{\operatorname{osc}_{#2}{#1}}

\usepackage{color}
\usepackage[colorlinks,citecolor=black,urlcolor=black]{hyperref}

\title[A scalar version of the Caflisch-Luke paradox]{A scalar version of the Caflisch-Luke paradox}

\author[A. Gloria]{Antoine Gloria}
\address[Antoine Gloria]{Sorbonne Universit\'e, CNRS, Universit\'e de Paris, Laboratoire Jacques-Louis Lions (LJLL), F-75005 Paris, France \& Universit\'e Libre de Bruxelles, Département de Mathématique, Brussels, Belgium}

\begin{document}
\selectlanguage{english}

\begin{abstract}
Consider an infinite cloud of hard spheres sedimenting in a Stokes flow in the whole space $\R^d$.
Despite many contributions in fluid mechanics and applied mathematics, there is so far no rigorous definition of the associated 
effective sedimentation velocity. Calculations by Caflisch and Luke in dimension $d=3$ suggest that the effective velocity is well-defined for
hard spheres distributed according to a weakly correlated and dilute point process, and that the variance of the sedimentation speed is infinite.
This constitutes the Caflisch-Luke paradox. In this contribution, we consider a scalar version of this problem that displays the same
difficulties in terms of interaction between the differential operator and the randomness, but is simpler in terms of PDE analysis.
For a class of hardcore point processes we rigorously prove that
 the effective velocity is well-defined in dimensions $d>2$, and that the variance is finite in dimensions $d>4$, confirming the 
formal calculations by Caflisch and Luke, and opening a way to the systematic study of such problems.
\end{abstract}

\maketitle

\section{Introduction and main results}

\subsection{The Caflisch-Luke paradox in sedimentation and its scalar version}

Sedimentation describes suspensions of rigid particles in a fluid that fall due to buoyancy (that is, particles are heavier than the fluid, and
fall by gravity), cf.~\cite{MR0014916,Batchelor}.
In physical experiments, heavy particles are spread in a tank, stirred, and then left falling at the bottom, cf.~\cite{sedimentation,pete}.
Assuming that the fluid is at rest at all times and that particles are identical balls, there are four parameters in the system 
at fixed time: the positions of the particles $\calP$, the size $r$ of the particles, the size $R$ of the tank, and the buoyancy $g$ (the 
effective gravity).
When the system ``reaches'' a stationary regime (the so-called constant composition zone) in the bulk of the tank, one can define an apparent effective sedimentation speed $\bar u_{\calP,r,R,g}$ as the average of the sedimentation speeds of all the particles.
One says that the sedimentation speed $\bar u_{\calP,r,g}$ is well-defined if $\bar u_{\calP,R,r,g}$ has a finite limit when $R\uparrow \infty$
(which, in experiments, means that $\bar u_{\calP,R,r,g}$ does not depend on the size $R$ of the tank when the tank 
is large enough), cf.~\cite{MR2768013}.
From a mathematical perspective this means one should be able to consider the thermodynamic limit of the problem, that is, consider
the limit $R\uparrow +\infty$ of the size of the tank, and define a limiting problem in the whole space $\R^d$.
In the following we write such a model in $\R^d$ and normalize most of the physical constants.
The only significant parameter left is the point set that describes the positions of particles. 
The effect of stirring the fluid is expected to reduce very much correlations between the particles at initial time.
Assumptions on the mixing properties on the point set will be crucial.

\medskip

Consider an infinite cloud of identical disjoint rigid spherical particles $B_i$ (the union of which we denote by $\calB=\cup_i B_i$) that sediment in a Stokes fluid in $\R^d$.
In the stationary regime, the velocities $u_i$ of the particles $B_i$ centered at $x_i$ (we set $\calP = \cup_i \{x_i\}$), the velocity $u$ of the fluid, and the pressure $p$
in the fluid satisfy
the coupled system of equations 
\begin{equation}\label{e.corrSed}
\left\{
\begin{array}{l}
-\triangle u\,=\, \nabla p \quad \text{ in } \R^d \setminus \calB,
\\
\nabla \cdot u\,=\,0 \quad \text{ in } \R^d \setminus \calB,
\\
\forall i:\quad u|_{B_i} \equiv u_{i} \in \R,  \quad \int_{\partial B_i} (\nablas u+p\Id) \cdot \nu =  g,
\end{array}
\right.
\end{equation}
where $g$ is the buoyancy (supposed constant) --- we have neglected the rotations of the particles.
The particles are ``active'' in the sense that they experience the gravity in a different way than the fluid (their density of mass is different).
If it exists, the effective sedimentation velocity $\bar u$ of the particles is given by
$$
\bar u = \lim_{R\to \infty} \frac{\sum_{i \in \calP \cap B^R} u_i}{|\calP \cap B^R|},
$$
where $B^R=B(0,R)$ is the ball of radius $R$ centered at the origin.
If $\calP$ is a stationary ergodic point set and \eqref{e.corrSed} is well-posed, one expects by stationarity and ergodicity $\bar u=\expec{u_i}$ for any $i$.
In fluid mechanics, determining the effective sedimentation velocity when the density of the particles is small is known as the Batchelor problem.
Despite several contributions in fluid mechanics \cite{MR0014916,Batchelor} and applied mathematics \cite{MR815929,MR1021642}, \emph{there is so far no proper definition of $\bar u$}.
What calculations by Caflisch and Luke \cite{Caflisch-Luke} suggest is that for point processes that are weakly correlated (in some sense),
and in the regime of low density $\theta \ll 1$ of particles,
\begin{itemize}
\item $\bar u$ is well-defined for $d=3$,
\item $\expec{u_i^2}=\infty$ for $d=3$.
\end{itemize}
These two (conjectured) properties constitute the Caflisch-Luke paradox: although the effective sedimentation velocity is well-defined, the associated variance is infinite. Understanding why this paradox is not observed in experiments remains an active field of research in fluid mechanics, cf.~\cite{MR2768013}.

\medskip

Let us now comment a bit on \eqref{e.corrSed} as a random PDE on $\R^d$. The differential operator is deterministic (it is the Stokes operator), the domain is random (it is $\R^d$ minus 
the union $\calB$ of particles), the boundary conditions on the particles depend nonlinearly and nonlocally on the point set $\calP$.
The difficulty in this equation is twofold: the map $\calP \mapsto \{u_i\}_i$ is nonlinear and nonlocal and the randomness appears in a lower-order term of the operator (it is at the level of $u_i$, not in the Laplacian). The first difficulty is reminiscent of the corrector equation in stochastic homogenization, which is by now well-understood (cf.~\cite{GO1,GNO1,Gloria-Otto-10b,GNO-quant},
and also \cite{AS,GO4,AKM2,AKMbook}) -- for instance using functional inequalities in probability as we shall do here. The second difficulty is reminiscent of the Schr\"odinger operator in a random potential. Whereas the PDE analysis is more involved in homogenization (the higher-order operator does not have constant coefficients), the difficulty is more on the probability side for random Schr\"odinger operators (there is less averaging in a lower-order term, cf.~Remark~\ref{rem:div-or-not} below).
A third difficulty is the incompressibility constraint (and therefore the pressure). 
From a probabilistic point of view, this difficulty is not essential:  the pressure is obtained by taking a Helmholtz projection,
which hardly amplifies correlations, and can therefore be neglected in front of the effect of the zero-order term --- which is 
why we consider this simpler model here. This additional difficulty for Stokes is therefore only on the PDE analysis side.
Yet, from the physical point of view, pressure allows to equilibrate forces, which means that if we want to neglect the pressure, we need to put a back flow into the picture.

\medskip

The aim of the present work is to investigate  the Caflisch-Luke paradox for a simpler equation that shares
the same basic difficulties as the sedimentation problem \eqref{e.corrSed} in terms of interaction of the differential operator with the randomness. 
We consider the following scalar equation posed on the whole space (which, in line with stochastic homogenization, we call corrector equation)
\begin{equation}\label{e.corr}
-\triangle u\,=\,\bar g_\theta \quad \text{ in } \R^d\setminus\calB, \quad  \forall i:\quad u|_{B_i} \equiv u_i \in \R,  \quad \int_{\partial B_i} \partial_n u = \bar g ,
\end{equation}
where $\bar g \in \R$ is given, $B_i$ are disjoint spherical inclusions centered at $x_i$ and of unit volume, $\calP= \cup_i \{x_i\}$ is a stationary ergodic point process,
$\theta = \expec{\mathds 1_{\calB}}$ is the intensity of the point set (equivalently, the density of inclusions), and with the neutrality condition $\bar g_\theta (1-\theta)- \bar g \theta=0$ (that is, $\bar g_\theta = \frac{\theta}{1-\theta}  \bar g$).
The unknown is the function $u \in H^1_\loc(\R^d)$ (and therefore the $u_i$'s).
In particular, under which conditions on the point set $\calP$, dimension $d$, and $\theta$ are
\begin{itemize}
\item[(I)] the ``corrector problem'' \eqref{e.corr} well-posed? 
\item[(II)] the effective electric field $\bar u=\expec{u_i}$ well-defined?
\item[(III)] the variance $\expec{u_i^2}$ of the electric field well-defined or infinite?
\end{itemize}

\medskip

The present approach towards sedimentation of particles considers the regime when particles very strongly interact, albeit in a stationary regime.
For results on dynamical aspects of sedimentation (either up to times for which particles do not strongly interact or in some homogenization regime), see \cite{JO-04,Hofer-18}.
For the related (but technically quite different) problem of justification of the effective viscosity due to ``passive'' particles in a Stokes flow (which is the case when the density of mass of the fluid and of the particles is the same), we refer the reader to the recent works \cite{HM-12,GVH-19}.

\medskip

Notation:
\begin{itemize}
\item For all (unit volume) inclusions $B_i$ centered at $x_i$ and all $t\ge 1$ we denote by $B_i^t:=B(x_i,t)$ the ball of radius $t$ centered at $x_i$;
\item $\lesssim$ (resp. $\gtrsim$) means $\le \times C$ (resp. $\ge \times C$) for some constant $C$ depending only on dimension (if not otherwise stated via a subscript on $\lesssim$). When both $\lesssim$ and $\gtrsim$ hold, we write $\sim$. When the multiplicative constant needs to be large enough, we write $\ll$ (resp. $\gg$).
\end{itemize}

\subsection{Massive approximation and main results} \label{sec:massive}

As standard in stochastic homogenization, we 
introduce a massive approximation of the corrector equation, and add an infra-red regularization which aims at localizing the problem.
Let $T\gg 1$, and consider on the whole space
\begin{equation}\label{e.corr-T}
\frac1Tu_T-\triangle u_T\,=\, \bar g_\theta \quad \text{ in } \R^d \setminus \calB, \quad  \forall i:\quad u_T|_{B_i} \equiv u_{T,i} \in \R,  \quad \int_{\partial B_i} \partial_n u_T = \bar g.
\end{equation}
Existence and uniqueness of solutions are proved on a deterministic basis in the following lemma,
as well as the finiteness of the massive effective electric field $\bar u_{T}:=\expec{u_T \mathds{1}_\calB}\theta^{-1}$.

\medskip

Before we state this result, let us recall some standard notions in PDEs with random coefficients.
Let $L^1(\Omega,\mathcal F,\mathbb P)$ denote a probability space, where $\Omega=\calS_1$ is the set of hardcore point sets $\calP$ (seen as an infinite sum of Dirac masses, endowed with the topology associated with the duality with continuous functions), distributed according to some stationary and ergodic probability measure $\mathbb P$, and $\expec{\cdot}$ the associated expectation.
Stationarity means that for all $y\in \R^d$, $\calP$ and $T_y \calP:=y+\calP=\{x+y\,|\,x\in \calP\}$ (which belongs to $\Omega$ by definition) have the same joint distribution under $\mathbb P$. Ergodicity means that if an event $E \in \mathcal F$ is such that $T_y E \subset E$ for all $y\in \R^d$, then $\pr{E} \in \{0,1\}$.
A random variable $Y$ is a measurable function on $\Omega$. A (jointly measurable) random field $Y \in L^1(\Omega,L^1_\loc(\R^d))$ is
said to be stationary if for all $x,y\in \R^d$, we have $Y(\cdot+y,x)=Y(\cdot,x+y)$ almost surely. Given $Y\in L^1(\Omega)$, the stationary extension $\tilde Y\in L^1(\Omega,L^1_\loc(\R^d))$ of $Y$ is defined by $\tilde Y(\calP,x):=Y(\calP+x)$. As customary in the field, we shall not distinguish between $Y$ and $\tilde Y$ (and use the same notation).
\begin{lem}\label{lem:massive-random}
Let $\rho \ge \sqrt{d}+1$, $\calP=\{x_i\}_{i\in \N}$ be a random stationary ergodic point set taking values in $\calS_\rho=\{ \calP' \in (\R^d)^\N \,:\, \forall x\ne x' \in \calP', |x-x'|\ge \rho\}$.
Let $\bar g \in \R$ and recall that $\theta =\expec{\mathds{1}_{\calB}}$.
For all $T>0$, there exists a unique stationary field $u_T$ that almost surely belongs to $\Hc_\uloc:=\{v\in H^1_\loc(\R^d),  \forall i:v|_{B_i}\equiv v(x_i),  \sup_{x\in \R^d} \int_{B(x)} v^2+|\nabla v|^2<\infty\}$  
and solves \eqref{e.corr-T} in the distributional sense, that is, for all stationary fields $v \in L^2(\Omega, \Hc_\uloc)$, 
\begin{equation}\label{e.eq-weak-prob}
\frac1T \expec{u_T v \mathds1_{\R^d\setminus\calB}}+\expec{\nabla u_T \cdot \nabla v} \,=\, \expec{v (\bar g_\theta \mathds1_{\R^d\setminus\calB}-\bar g \mathds1_{\calB})}.
\end{equation}
In particular, 
\begin{equation}\label{e.eq-expec-uT}
 \expec{u_T \mathds1_{\R^d\setminus\calB}} \,=\,  0
\end{equation}
and $u_T$ satisfies the energy estimate
\begin{equation}\label{e.massive-apriori-random}
\expec{\frac1T u_T^2 +  |\nabla u_T|^2 } \,\lesssim \, T  \theta \bar g^2
\end{equation}
and the identity
\begin{equation}\label{e.expec-V_iT}
\expec{(\frac1T u_T^2+|\nabla u_T|^2)\mathds1_{\R^d\setminus\calB}} = -\bar g\expec{u_T\mathds 1_{\calB}}. 
\end{equation}
\end{lem}
The next natural step is to pass to the limit  $T\uparrow \infty$ in \eqref{e.corr-T} to recover \eqref{e.corr}.
The energy estimate \eqref{e.massive-apriori-random} is however not enough --- not even for $\nabla u_T$, which contrasts very much with the corrector 
equation in stochastic homogenization for which existence of stationary gradients comes for free.
This lack of compactness is the main issue in the analysis of~\eqref{e.corr}.
Compactness will be obtained the hard way using a quantitative assumption of ergodicity and several regularity estimates.

\medskip

Due to the nonlinear dependence of $u_T$ with respect to $\calP$, it is natural to assume quantitative ergodicity in form of a functional inequality in probability. Indeed, functional inequalities in probability provide a calculus which allows one to linearize the dependence on the randomness.
In view of the hardcore condition, standard functional inequalities such as spectral gap do not apply, and we need to appeal to the multiscale functional inequalities in probability
introduced by Duerinckx and the author in \cite{DG1}.
\begin{defin}[Multiscale Poincar\'e inequality \cite{DG1}]\label{def:SG}
Let $\rho \ge 1$. We say that a point process $\calP$ taking values in $\calS_\rho$ 
satisfies a multiscale Poincar\'e inequality in probability if there exists $C<\infty$ such that for all measurable functions $Y:\calS_\rho\to \R$ we have 
\begin{equation}\label{e.SG}
\var{Y} \,\le\, C  \int_1^\infty \int_{\R^d} \expec{\oscd{Y}{B^\ell(z)}}dz e^{-\frac \ell {C}}  d\ell,
\end{equation}
where $Y$ is understood as the random variable $Y(\calP)$, and the oscillation on some subset $D\subset \R^d$ 
is defined by
$$
\osc{Y}{D}(\calP_\rho)\,:=\sup_{\calP'\in \calS_\rho \,:\, \calP'|_{\R^d\setminus D}\equiv \calP|_{\R^d\setminus D}} Y(\calP')-Y(\calP).
$$
\end{defin}
Two typical examples of point processes satisfying \eqref{e.SG} are the hardcore Poisson point processes and the random parking measure (both defined 
for $\rho\ge 1$ via the Penrose graphical construction \cite{Penrose-01} starting from the Poisson point process of intensity unity on $\R^d\times [0,\lambda]$ and on $\R^d\times \R_+$, respectively, for some $\lambda >0$), cf.~\cite{DG1}. 
Let $J$ denote the jamming limit defined in \cite{Penrose-01} (that is, the density of spherical inclusions of radius $1$ centered at points of the random parking measure
of parameter $\rho=1$), and $B(0)$ be the ball of unit volume centered at 0. We recall the three main geometric properties of these points sets:
\begin{itemize}
\item Hardcore Poisson process $\calP$ of parameters $(\rho,\lambda)$: 
\begin{multline*}
\inf\{\dist(x,\calP\setminus \{x\})\,:\,x\in \calP\}\ge \rho, \quad \sup\{\dist(x,\calP\setminus \{x\})\,:\,x\in \calP\}=\infty, 
\\
\expec{\sum_{x\in \calP} \mathds{1}_{B(0)}(x)} \,\le\, \lambda \wedge (J\rho^{-d});
\end{multline*}
\item Random parking measure $\calP$ of parameter $\rho$:
\begin{multline*}
\inf\{\dist(x,\calP\setminus \{x\})\,:\,x\in \calP\}\ge \rho, \quad \sup\{\dist(x,\calP\setminus \{x\})\,:\,x\in \calP\}\lesssim \rho,
\\
 \expec{\sum_{x\in \calP} \mathds{1}_{B(0)}(x)} \,=\, J\rho^{-d}.
\end{multline*}
\end{itemize}

\medskip

The main result of this article is the existence of solutions for \eqref{e.corr}:
\begin{theor}\label{th:effective}
Let $\bar g =1$ and $d>2$.
There exist  $\rho_\m \ge \sqrt d+1$ such that 
if  $\calP$ is a hardcore point process of parameter $\rho \ge \rho_\m$ that satisfies  \eqref{e.SG}, then \eqref{e.corr} admits a unique solution $u$ the gradient of which is stationary and has finite second moment, and the effective electric field $\bar u =\lim_{T\uparrow +\infty} \bar u_{T}$ 
 is well-defined.
In addition, the solution $u_T$ of \eqref{e.corr-T} satisfies $\lim_{T\to \infty}\expec{|\nabla u_T-\nabla u|^2}=0$ and $\lim_{T\to \infty} \frac1T \expec{u_T^2}=0$.
\end{theor}
This result is completed by the following Caflisch-Luke estimates. 
\begin{prop}\label{prop:unif-uT}
Under the assumptions of Theorem~\ref{th:effective} the solution $u_T$ of \eqref{e.corr-T} satisfies 
for all $p\ge 1$, $\sqrt T\gg \rho$, and $d>2$ the estimates 
\begin{equation}\label{e.unif-uT-0}
\expec{ \Big(\fint_{B} |\nabla u_T|^2\Big)^p}^\frac{1}{p} \, \lesssim_\rho \, p^{\frac{d+2}{d-2}\gamma},
\end{equation}
and
\begin{equation}\label{e.unif-uT-1}
\expec{ \Big(\fint_{B} u_T^2\Big)^p}^\frac{1}{p} \, \lesssim_\rho \,  \left\{
\begin{array}{rcl}
d=3&:&p^{\frac{7}3  \gamma} \sqrt T,
\\
d=4&:&p^{2^+\gamma }  \log  T,
\\
d>4&:&p^{\frac{d+4}{4}\gamma },
\end{array}
\right.
\end{equation}
where  $\gamma:=d+5+\frac {(d-1)(d+2)}{d}$, and $2^+$ means any real number larger than $2$.
\end{prop}
Theorem~\ref{th:effective} and Proposition~\ref{prop:unif-uT} solve questions (I)--(III).
In particular, if the point set $\calP$ satisfies a suitable functional inequality, we have \emph{with a largeness condition} on the hardcore parameter $\rho$  but \emph{without (additional) smallness condition} on the intensity of $\calP$ that
\begin{itemize}
\item \eqref{e.corr}
is well-posed and the effective electric field $\expec{u_i}$ is well-defined for $d>2$;
\item $\expec{u_i^2}$ is finite for $d>4$.
\end{itemize}
Proposition~\ref{prop:unif-uT} suggests that $\expec{u_T^2}$ does not
remain uniformly bounded wrt $T$ in dimensions 3 and 4.
This supports the Caflisch-Luke paradox for $d=3$ and $d=4$, and rigorously shows that there is no paradox in dimensions $d>4$.
The restriction on $\rho$ in these results entails a deterministic positive distance between particles.
It is related to the deterministic regularity estimates of Lemmas~\ref{lem:reg1} and~\ref{lem:Green} below. 
We believe this condition might be relaxed provided one develops a random large-scale version of  Lemmas~\ref{lem:reg1} and~\ref{lem:Green} in the spirit of \cite{AS,GNO-reg}
in homogenization.

\subsection{Proposition~\ref{prop:unif-uT}, stochastic homogenization, and Coulomb interactions}\label{rem:div-or-not}

The estimates of Proposition~\ref{prop:unif-uT} are in line with the intuition that $u_T$ essentially behaves like the solution $v_T$ of
\begin{equation}\label{e.model1}
(\frac1T-\triangle) v_T = \mathds{1_\calB}-\theta \ \mbox{ in }\R^d,
\end{equation}
which can be seen as the linearization of \eqref{e.corr-T} with respect to the randomness $\calP$, and for which the corresponding estimates follow from an explicit calculation using the massive Green's function  $G_T$ (in this case these estimates are also lower bounds):
\begin{eqnarray*} 
\expec{v_T^2} &= &\int \int G_T(x)G_T(x') \cov{1_\calB(x)}{1_\calB(x')}dxdx',
\\
\expec{|\nabla v_T|^2} &= &\int \int \nabla G_T(x)\cdot \nabla G_T(x') \cov{1_\calB(x)}{1_\calB(x')}dxdx',
\end{eqnarray*}
so that \eqref{e.unif-uT-0} and \eqref{e.unif-uT-1} follow from the two-sided estimates 
\begin{eqnarray*}
d>2&:&|x|^{2-d}\exp(-\frac{|x|}{C})\lesssim G_T(x)\lesssim |x|^{2-d}\exp(-\frac{|x|}{c}),
\\
d\ge 2&:& |x|^{1-d}\exp(-\frac{|x|}{C})\lesssim |\nabla G_T(x)| \lesssim |x|^{1-d}\exp(-\frac{|x|}{c})
\end{eqnarray*}
provided $ \cov{1_\calB(x)}{1_\calB(x')}$ decays fast enough as $|x-x'|\to \infty$.
This yields in particular the sharp bounds for the Poisson point process 
\begin{equation}\label{e:easy-cal}
\expec{|\nabla v_T|^2} \,\simeq \, 
 \left\{
\begin{array}{rcl}
d=2&:&\log^\frac12 T,
\\
d>2&:&1.
\end{array}
\right.
\end{equation}
Note that the corrector equation in stochastic homogenization \cite[Appendix]{GO1} rather  behaves like the solution $w_T$ (for some unit vector $e\in \R^d$) of 
\begin{equation}\label{e.model2}
(\frac1T-\triangle) w_T = \nabla \cdot (\mathds{1_\calB}e) \ \mbox{ in }\R^d,
\end{equation}
for which we have \cite{Gloria-Otto-10b} for all $d\ge 1$
\begin{equation}\label{e.unif-uT-hom}
\expec{ \Big(\fint_{B} |\nabla w_T|^2\Big)^p}^\frac{1}{2p} \, \lesssim_p \, 1, \quad
\expec{ \Big(\fint_{B} w_T^2\Big)^p}^\frac{1}{2p} \, \lesssim_p \, \left\{
\begin{array}{rcl}
d=1&:& \sqrt T,
\\
d=2&:&\log^\frac12 T,
\\
d>2&:&1.
\end{array}
\right.
\end{equation}
The difference of scalings between \eqref{e.unif-uT-0}--\eqref{e.unif-uT-1} \& \eqref{e.unif-uT-hom} comes from the fact
that the RHS in \eqref{e.model2} \& \eqref{e.model1} is in divergence form or not.

\medskip

Equations like \eqref{e.model1} have already been considered in the setting of Gibbs measures of particles interacting
via the Coulomb potential (or more general Riesz potentials), see~\cite{Serfaty-14} for lecture notes by Serfaty on the context. In this case, $\mathds{1_\calB}$ is replaced by a sum of Dirac masses (and $\mathds{1_\calB}$ can be interpreted as a local smoothing), and the massive term by another screening procedure.
Provided one can pass to the limit $T\uparrow +\infty$ and define $\nabla v=\lim_{T\uparrow +\infty} \nabla v_T \in L^2(\Omega)$
(which follows from \eqref{e.unif-uT-1} for $d>2$, cf.~the proof of Theorem~\ref{th:effective}), one can define the energy of the (smeared out) point process $\mathds{1_\calB}$ as
$$
W_{\elec}(\mathds{1_\calB})\,:=\, \expec{|\nabla v|^2} \,=\, \lim_{T\uparrow +\infty} \expec{\frac1T v_T^2+|\nabla v_T|^2}.
$$
In particular, up to the singularity of the Dirac masses (cf.~the renormalized energy of Sandier and Serfaty \cite{Sandier-Serfaty-15}), 
\eqref{e:easy-cal} implies that the energy $W_{\elec}$ of the Poisson point process is infinite in dimension $d=2$ and finite in dimensions $d>2$,  as was also proved in~\cite{Leble}. 
In~\cite{Serfaty-14} and subsequent works, one of the primary goals is to minimize the functional $W_{\elec}$ (or a variant of it at finite temperature) on the set of stationary point sets and to obtain explicit formulas.
In the present contribution, the point of view is different: We give ourselves a stationary point set and prove the finiteness of the energy
 $\expec{(\frac1T u_T^2+|\nabla u_T|^2)\mathds1_{\R^d\setminus\calB}}$, where $u_T$ displays a nonlinear dependence wrt $\calB$.
 We do not address the minimization of this energy on stationary point sets, which we could reformulate as 
finding the stationary point set of fixed intensity that minimizes the electric field $\expec{u_i}=\lim_{T\uparrow +\infty}\expec{(\frac1T u_T^2+|\nabla u_T|^2)\mathds1_{\R^d\setminus\calB}}$ ---  which one might conjecture to be a crystal,
cf.~\cite{Theil-06,Blanc-Lewin-15,Petrache-Serfaty-19} for related models in dimensions $d\ge 2$.
Note that in the work \cite{GVH-19}, an energy and similar techniques as in~\cite{Serfaty-14} have been used in the context
of the Einstein formula for dilute suspensions in a Stokes fluid.

\bigskip

The rest of the paper is organized as follows.
In Section~\ref{sec:structure}, we display the structure of the proof, which is partly inspired by \cite{Gloria-Otto-10b}.
The key result  is the decay of averages of $(\frac1{\sqrt T} u_T,\nabla u_T)$ provided by Proposition~\ref{prop:unif-nablauT}, which allows to buckle and pass to the limit as $T\uparrow \infty$.
The proof of this result (which is displayed in Section~\ref{sec:main-results}, together with the proofs of Theorem~\ref{th:effective} and Proposition~\ref{prop:unif-uT})
follows from the combination of deterministic results (energy estimates, a compactness result, regularity results)
with the multiscale Poincar\'e inequality through a sensitivity calculus.
The regularity results are proved in Section~\ref{sec:regularity-results}, whereas the other auxiliary results are proved in Section~\ref{sec:auxiliary-results}.

\tableofcontents

%%%%%%%%%%%%%%%%%%%%%%%%
%%%%%%%%%%%%%%%%%%%%%%%%
%%%%%%%%%%%%%%%%%%%%%%%%

\section{Structure of the proofs}\label{sec:structure}

In the rest of this article, we assume that $\rho \ge \sqrt d+1$ (so that particles are at least at distance~1 from one another), that $T\gg \rho$, and we consider wlog $\bar g=1$.

\medskip

We start with the well-posedness of the massive approximation of the corrector equation,
in form of the following deterministic result.
\begin{lem}\label{lem:massive}
Let $g_1,g_2 \in \R$.
For all $T>0$ and all points sets $\calP =\{x_i\}_i \in \calS_\rho$, there exists a unique distributional solution $v_T$ of
\begin{equation}\label{e.corr-Tdet}
\frac1Tv_T-\triangle v_T\,=\, g_2 \quad \text{ in } \R^d \setminus \calB, \quad  \forall i:\quad v_T|_{B_i} \equiv v_{T,i} \in \R,  \quad \int_{\partial B_i} \partial_n v_T = g_1.
\end{equation}
in $\Hc_\uloc$. It satisfies the energy estimate
\begin{equation}\label{e.massive-apriori}
\sup_{x \in \R^d} \Big(\frac1T \fint_{B^{\sqrt T}(x) } v_T^2 +\fint_{B^{\sqrt T}(x) } |\nabla v_T|^2 \Big) \,\lesssim \, T  (g_1^2 +g_2^2).
\end{equation}
\end{lem}
Applied to a random stationary ergodic point process, it entails Lemma~\ref{lem:massive-random}, and therefore the existence and uniqueness
of the stationary field $u_T$.
In order to prove estimates on $u_T$ and $\nabla u_T$ that are uniform wrt $T$, it is natural to consider the random variable $Y=\Big(\int_{B(0)} |\nabla u_T|^2\Big)^\frac12$
and to apply the variance estimate \eqref{e.SG}.
This random variable is unfortunately not linear wrt $\nabla u_T$, which prevents us from using efficiently the linearity of the PDE \eqref{e.corr-T} to estimate differences of solutions
(as required by the oscillation in the RHS of \eqref{e.SG}).
The following lemma shows however that it is enough to apply the variance estimate to quantities of the form 
$Y=\int_{\R^d} \nabla u_T \cdot g$ for a finite number of (deterministic and) compactly supported functions $g$. 
This can be seen as a
compactness result for solutions of  \eqref{e.corr-T}, in the spirit of \cite{Gloria-Otto-10b}.
\begin{lem}\label{lem:compactness}
For all $\delta>0$, there exist a finite family $\{g_n\}_{1\le n\le N}$ (with $N$ depending on $\delta$)
of bounded vector-valued functions supported in $B^{2}$ normalized in $L^2(\R^d)^d$ and a constant $C<\infty$ such that we have
for all $R\ge 4 \rho$ and all $T>0$
\begin{multline}\label{e.compact0}
\int_{B^R} |\nabla u_T|^2 \,\le \, C \sum_{n=1}^N \Big(\int_{B^{2R}}  \nabla u_T\cdot R^{-\frac d2} g_n(\tfrac \cdot R)\Big)^2 + \delta \int_{B^{2R}} |\nabla u_T|^2 
\\
+C R^2 \int_{B^{2R}} (\frac{1}{T^2}u_T^2+ 1).
\end{multline}
In particular for all $p\ge 1$,  and all $\sqrt T\gg R$,
\begin{multline}\label{e.compact}
\expec{\Big(\fint_{B^R} |\nabla u_T|^2\Big)^p}^\frac1p \,\le \, C  \sum_{n=1}^N \expec{\Big|\int_{B^{2R}}  \nabla u_T\cdot g_{R,n}\Big|^{2p} }^\frac1p
\\
+C  R^2 
+ C \expec{\Big|\fint_{B^{2R}} \frac1{\sqrt T} u_T\Big|^{2p}}^\frac1p,
\end{multline}
with the short-hand notation $g_{R,n}:x \mapsto R^{-d} g_n(\tfrac x R)$.
\end{lem}
\begin{rem}
Although we can pass to the limit in the first RHS term of \eqref{e.compact} by the Birkhoff ergodic theorem,
the second RHS term blows up as $R\uparrow \infty$ so that  \eqref{e.compact} does not
yield the uniform boundedness of $\expec{|\nabla u_T|^2}$ without further information.
\end{rem}
The four upcoming results are solely based on PDE analysis, and will be combined to the multiscale Poincar\'e inequality~\eqref{e.SG} in order to prove 
Proposition~\ref{prop:unif-nablauT} (see below).

\medskip
For all $D \subset \R^d$, $\calP,\calP'\in \calS_\rho$ such that $\calP'|_{\R^d\setminus D}\equiv \calP|_{\R^d\setminus D}$,
we use the short hand notation $Y=Y(\calP)$, $Y'=Y(\calP')$, and $\delta_D Y=Y'-Y$.
We call $\{B_i\}_i$ and $\{B_i'\}_i$ the inclusions associated with $\calP$ and $\calP'$, respectively, and define 
$\calB,\calB'$ and $\R^d\setminus{\calB},\R^d\setminus{\calB'}$ accordingly.
Let $u_T$ be the solution of \eqref{e.corr-T} associated with $\calP$, $u_T'$ be the solution associated with $\calP'$, and set $w_T:=u_T'-u_T$.
The following lemma establishes the equation satisfied by $w_T$, as well as some energy estimate.
\begin{lem}\label{lem:diff-u}
The map $w_T$ solves on $\R^d$ the equation 
\begin{equation}\label{e.diff-uT}
\frac1T w_T-\triangle w_T\,=\, \bar g_\theta (\mathds 1_{\R^d\setminus \calB'}-\mathds 1_{\R^d\setminus \calB})+\frac1T(u'_T \mathds1_{\calB'}-u_T \mathds1_{\calB})
-\sum_i (\nu_i' -\nu_i ),
\end{equation}
where $\nu_i= \nabla u_T \cdot n_i \delta_{\partial B_i}$ (resp. $\nu_i'=\nabla u_T' \cdot n'_i \delta_{\partial B_i'}$),
with $n_i$ (resp. $n'_i$) the outward normal of $B_i$ (resp. $B_i'$), are measures that belong to $H^{-1}(\R^d)$, and satisfy  satisfy $\nu_i(1)=\nu_i'(1)=1=|\partial B|\tilde g$ for all $i$.
For $d>2$, $w_T$ satisfies the estimate
\begin{multline}\label{e.diff-uT-estim}
\int_{\R^d}\frac1T w_T^2+|\nabla w_T|^2 \,\le\,  C |\calB\triangle \calB'|^\frac{d+2}{d}+C  |\calB\triangle \calB'|^{\frac 2d} \frac1{T^2} \int_{\calB\triangle \calB'} u_T^2+C \int_{\overline D} |\nabla u_T|^2
\\
+\sum_{B_i \in \calB \setminus  \calB'} C\Big(\int_{\partial B_i} |\nabla u_T \cdot n_i-\tilde g|\Big)^2
+\sum_{B_i' \in \calB' \setminus \calB} C\Big(\int_{\partial B_i'} |\nabla u_T' \cdot n_i'-\tilde g|\Big)^2
\end{multline}
for some $C\gg 1$, with the notation 
$\calB \triangle \calB'=\{x \in \calB\,|\, x \notin \calB'\} \cup \{x\in \calB'\,|\, x\notin \calB\}$, 
 the less standard notation $\calB\setminus \calB':=\cup_{B\subset \calB, B\not \subset \calB'}B$ (in particular, $\calB\setminus \calB'$ is a union of balls that might intersect  but 
are not included in $\calB'$) and $\calB'\setminus \calB:=\cup_{B'\subset \calB', B'\not \subset \calB}B'$,
and where $\overline{D}$ is a short-hand notation for the enlarged set $\{x \in \R^d\,|\, \dist(x,D)\le \sqrt d+1\}$.
\end{lem}
In order to use \eqref{e.diff-uT-estim} to control the RHS of \eqref{e.SG}, we need two regularity results.
The first result is the following quantitative estimate of approximate radiality.
\begin{lem}\label{lem:reg1}
For $d>2$ there exists an exponent $\alpha>0$ (coming from hole-filling)
such that for all $\rho\ge 4$, $T\ge 1$, all $g_1,g_2,\bar v \in \R$, and all $v \in H^1(B^\rho)$ that satisfy
$$
\frac1T v-\triangle v\,=\, g_2 \text{ in }B^\rho\setminus B^1, \quad v\equiv \bar v \text{ on }\partial B^1, \quad \fint_{\partial B^1} \nabla v\cdot n=g_1,
$$
where $B^1$ and $B^\rho$ denote the balls of radius $1$ and $\rho$, respectively (both centered at the origin), we have  
$$
\int_{\partial B^1} \Big|\nabla v\cdot n-g_1\Big|\,\lesssim\, \rho^{-\alpha} \Big(\int_{B^\rho} |\nabla v|^2 +\rho^{d}(1+\frac{\rho^2}{T})(g_1^2+ \rho^2 (g_2^2+\frac{\bar v^2}{T^2}) ) \Big)^\frac12.
$$
\end{lem}
The combination of Lemmas~\ref{lem:diff-u} and~\ref{lem:reg1} completes the a priori estimate for $w_T$
in terms of $u_T$ only.
\begin{cor}\label{cor:apriori}
For all $d>2$ there exists $\rho_\m\ge \sqrt d+1$ such that for all $\rho>\rho_\m$ (cf.~$\calS_\rho$) and $T\gg \rho^{d+4}$, the map $w_T$ satisfies the estimate
\begin{equation}\label{e.diff-uT2}
\int_{\R^d}\frac1T w_T^2+|\nabla w_T|^2 \,\lesssim_\rho \,     |\calB \triangle \calB'|^\frac{d+2}{d}+  |\calB\triangle \calB'|^{\frac 2d} \frac1{T^2} \int_{\overline D} u_T^2
+ \int_{\overline{D}} |\nabla u_T|^2 ,
\end{equation}
where $\overline{D}$ is now a short-hand notation for the enlarged set $\{x \in \R^d\,|\, \dist(x,D)\le \rho_\m\}$.
%,
\end{cor}
The second regularity result we need is the following decay of Green's functions and of their first gradients.
\begin{lem}\label{lem:Green}
For $d>2$, there exists $\rho_\m\ge \sqrt d+1$ depending only on $d$ such that for all $\rho>\rho_\m$, all $T \gg \rho^2$, and $\calP \in \calS_\rho$,
the Green's function $x\mapsto G_T(x,y)$ defined for all $y \in \R^d\setminus \cup_{i} B^2_i$ as the unique
solution of
\begin{equation}\label{e.corr-Green}
\frac1T G_T(x,y)-\triangle G_T(x,y)\,=\, \delta(y-x)\ \text{ in } \R^d  \setminus \calB, \,  \forall B \subset \calB:\, G_T|_{B} \equiv G_{T,B} \in \R,  \ \int_{\partial B} \partial_n G_T = 0
\end{equation}
satisfies the pointwise estimates for all $x\in \R^d \setminus B^{2}(y)$
\begin{eqnarray}
|G_T(x,y)|&\lesssim & |x-y|^{2-d} \exp(-\frac{|x-y|}{C_d \sqrt T}),\label{e.Green1}
\\
|\nabla_x G_T(x,y)|&\lesssim &|x-y|^{1-d} \exp(-\frac{|x-y|}{C_d \sqrt T}),\label{e.Green2}
\end{eqnarray}
where the multiplicative constant only depends on $d$ (through $\rho_\m$).
\end{lem}
From now on, we call $\rho_\m$ the largest of the two radii defined in Lemmas~\ref{lem:reg1} and~\ref{lem:Green}.
Theorem~\ref{th:effective} and Proposition~\ref{prop:unif-uT} will follow from the upcoming result, that allows us to buckle in \eqref{e.compact}.
\begin{prop}\label{prop:unif-nablauT}
Let  $d>2$.
Let $g \in L^\infty(\R^d)^{d+1}$ be supported in $B^2$, and for all $R>0$ set $g_R:=R^{-d}g(\tfrac \cdot R)$.
Then for all $\rho \ge \rho_\m$, for all $p\ge 1$, $T \gg \rho^{d+4}$, and all $R\gg 1$, we have with $\gamma=d+5+\frac {(d-1)(d+2)}{d}$
\begin{multline}\label{e.unif-nablauT}
\expec{\Big|\int_{B^{2R}}  (\frac1{\sqrt T}u_T,\nabla u_T)\cdot g_{R}\Big|^{2p} }^\frac1p \\
\lesssim_\rho \, p^{\gamma} R^{2-d} \Big(1
+\frac1{T} \expec{\Big|\fint_{B^{2R}} \frac1{\sqrt T}u_T\Big|^{2p}}^\frac1p+ (1+\frac{R^2}{T^2})\expec{\Big(\fint_{B^{R}}  |\nabla u_T|^2\Big)^{p} }^\frac1p\Big).
\end{multline}
\end{prop}
In order to obtain the optimal power of the logarithm in the Caflisch-Luke estimate in dimension $d=4$,
we need the following slight refinement of Proposition~\ref{prop:unif-nablauT}.
\begin{cor}\label{cor:unif-nablauT}
Let $d>2$ and for all $R\gg 1$ let $g_R:x\mapsto (1+|x|)^{1-d}\mathds1_{B^R}$.
Then for all $\rho \ge \rho_\m$, for all $p\ge 1$, $T \gg \rho^{d+4}$, and all $R\gg 1$, we have with $\gamma=d+5+\frac {(d-1)(d+2)}{d}$
\begin{equation}\label{e.unif-nablauT-crit}
\expec{\Big|\int  \nabla u_T \cdot g_{R}\Big|^{2p} }^\frac1p \,
\lesssim_\rho \, p^{\gamma} \mu_d(R)\expec{1+\Big(\frac1{T^2}\fint_{B^{2}}u_T^2\Big)^p+\Big(\fint_{B^{2}} |\nabla u_T|^2\Big)^p}^\frac1p,
\end{equation}
where 
$$
\mu_d(R)= \left\{
\begin{array}{rcl}
d=3&:& R,
\\
d=4&:&   \log  R,
\\
d>4&:&1.
\end{array}
\right.
$$
\end{cor}

\section{Proofs of the main results}\label{sec:main-results}

\subsection{Theorem~\ref{th:effective}: Existence and uniqueness of correctors}
We split the proof into three steps. We first prove existence and uniqueness of solutions of \eqref{e.corr} by approximation
with a massive term based on Proposition~\ref{prop:unif-uT}.
Then we prove the convergence of $\bar u_{T}$ to $-\theta^{-1}\expec{|\nabla u|^2}$, and finally address the 
strong convergence of $\nabla u_T$ to $\nabla u$.

\medskip

\step1 Existence and uniqueness of $u$.

By Proposition~\ref{prop:unif-uT}, since $\nabla u_T$ is bounded in $L^2(\Omega)^d$, there exists $u \in L^2_\loc(\R^d,L^2(\Omega))$
such that $\nabla u$ is stationary, has finite second moment, and such that  $\nabla u_T$ converges  
weakly in $L^2(\Omega)^d$ to $\nabla u$ along some subsequence (which we do not relabel).
By the bound \eqref{e.unif-uT-1} on $u_T$, we may pass to the limit in the weak formulation \eqref{e.eq-weak-prob}, which shows that 
for all stationary fields $v \in L^2(\Omega, \Hc_\uloc)$, we have
\begin{equation}\label{e.weak-form-u}
\expec{\nabla u \cdot \nabla v} \,=\, \expec{v (\bar g_\theta \mathds1_{\R^d\setminus\calB}- \mathds1_{\calB})}.
\end{equation}
Let $\tilde u$ be another solution of \eqref{e.corr} such that $\nabla \tilde u$ is stationary and has finite second moment.
By the same argument as for the proof of \eqref{e.eq-weak-prob}, $\tilde u$ also satisfies \eqref{e.weak-form-u} (cf.~proof of Lemma~\ref{lem:massive-random} below),
so that the difference $w=u-\tilde u$ satisfies for all stationary fields $v \in L^2(\Omega, \Hc_\uloc)$
\begin{equation}\label{must-be-zero}
\expec{\nabla w \cdot \nabla v}=0.
\end{equation}
Let us prove that necessarily $\nabla w\equiv 0$. 
For all $\mu>0$, let $\mathfrak M^\mu:C^\infty_c(\R^d)\to C^\infty_c(\R^d)$ be a map that 
modifies smooth functions on a $\mu$-neighborhood of $B(0)$ to make them constant in $B(0)$,
and for all $\calP \in \calS_\rho$ define $\calM^\mu_\calP:C_c^\infty(\R^d)\to C_c^\infty(\R^d)$ as
$\calM^\mu_\calP =\prod_i \mathfrak M^\mu(x_i+\cdot)$.
Let now $\chi \in C^\infty_c(\R^d)$, $\zeta \in L^\infty(\Omega)$ (that is, $\zeta$ is a function of point sets),  and 
 define $v:(\calP,x)\,\mapsto\, \int_{\R^d} (\calM^\mu_{y+\calP} \chi)(y-x) \zeta(y+\calP)dy$ (which is finite since $\chi$ has compact support). 
By construction, $v$ is stationary and belongs to $L^2(\Omega, \Hc_\uloc)$. Indeed, for all $z\in \R^d$,
\begin{eqnarray*}
v(\calP,x+z)&=&\int_{\R^d} \calM^\mu_{y+\calP} \chi(y-x-z) \zeta(y+\calP)dy \\
&\stackrel{y'=y-z}=& \int_{\R^d} \calM^\mu_{y'+(z+\calP)} \chi(y'-x) \zeta(y'+(z+\calP))dy'\,=\,v(z+\calP,x),
\end{eqnarray*}
and we have  
$$
\expec{{\sup_{x}}^2 \{|v(\calP,x)|+|\nabla v(\calP,x)|\}}^\frac12 \,\lesssim \,\mu^{-1} \|\chi\|_{W^{1,\infty}(\R^d)} |\supp \chi| \|\zeta\|_{L^\infty(\Omega)}
$$
so that $v\in L^2(\Omega, \Hc_\uloc)$.
Note that  $\calM^\mu$ is weakly continuous in the following sense: for all bounded domains $D$ and all sequences $(\chi_n)_n$ of functions compactly supported in $D$, if  $\chi_n \rightharpoonup \chi$ weakly in $H^1(D)$, then almost surely $\calM^\mu \chi_n \rightharpoonup \calM^\mu \chi$ in $H^1(D)$.
 We then use \eqref{must-be-zero} with this choice of $v$ and obtain by construction, stationarity of $\nabla w$,
and the stationarity of the probability measure
\begin{eqnarray*}
0&=&\expec{\nabla w \cdot \nabla v(0)}=\expec{\int_{\R^d} \nabla (\calM^\mu_{y+\calP} \chi(y)) \cdot \nabla w(-y,y+\calP) \zeta(y+\calP)dy}
\\
&=&\expec{\zeta(\calP)\int_{\R^d} \nabla (\calM^\mu_{\calP} \chi(y)) \cdot \nabla w(-y,\calP) dy}.
\end{eqnarray*}
By arbitrariness of $\zeta$ (and the density of $L^\infty(\Omega)$ in $L^2(\Omega)$) and 
the weak continuity of $\calM^\mu$, this implies that almost surely 
we have for all $\chi$ 
$$
\int_{\R^d} \nabla (\calM^\mu_{\calP} \chi(-y)) \cdot \nabla w(y,\calP) dy\,=\,0.
$$
By the arbitrariness of $\mu>0$ and of $\chi$, this implies that $\nabla w\equiv 0$ on $\R^d\setminus \calB$ almost surely,
whereas $\nabla w \equiv 0$ on $\calB$ since $w$ is constant on the inclusions.
Uniqueness is proved.

\medskip

\step2 Existence of $\bar u=\lim_{T\to \infty} \bar u_{T}$.

Let $T,T'\ge 1$.
The starting point is \eqref{e.eq-weak-prob} for $u_T$ and $v=u_{T'}$ combined with  \eqref{e.eq-expec-uT} in the form $\expec{u_{T'} \mathds1_{\R^d\setminus\calB}}=0$
to the effect that
\begin{equation}\label{e.weak-diag}
\frac1T \expec{u_T u_{T'} \mathds1_{\R^d\setminus\calB}}+\expec{\nabla u_T \cdot \nabla u_{T'}} \,=\, -\expec{u_{T'} \bar g \mathds1_{\calB}}=-\bar u_{T'} { \theta}.
\end{equation}
We first pass to the limit $T\uparrow +\infty$, which yields
$$
\expec{\nabla u \cdot \nabla u_{T'}} \,=\,  -\bar u_{T'}   { \theta}.
$$
Taking then the limit $T'\uparrow +\infty$ finally shows
\begin{equation}\label{e.def-Utheta}
\bar u:= \lim_{T'\to \infty} \bar u_{T'} = \lim_{T'\to \infty}  \tfrac1{ \theta}  \expec{\nabla u \cdot \nabla u_{T'}} =-\tfrac1{ \theta} \expec{|\nabla u|^2}.
\end{equation}

\medskip

\step3 Strong convergence of $\frac1T u_T$ and $\nabla u_T$.

On the one hand, by the weak lower-semicontinuity of the norm, we have
$$
\expec{|\nabla u|^2} \,\le \, \liminf_{T\to \infty} \expec{|\nabla u_T|^2}.
$$
On the other hand, by \eqref{e.def-Utheta} and \eqref{e.weak-diag} for $T'=T$,
$$
\lim_{T\to \infty} \frac1T \expec{u_T^2 \mathds1_{\R^d\setminus\calB}}+\expec{|\nabla u_T|^2} = \expec{|\nabla u|^2} .
$$
The combination of these two properties then implies
$$
\lim_{T\to \infty} \frac1T \expec{u_T^2 \mathds1_{\R^d\setminus\calB}}=0, \quad \lim_{T\to \infty} \expec{|\nabla u_T|^2} = \expec{|\nabla u|^2},
$$
which in turn yields the strong $L^2(\Omega)$ convergence of $\nabla u_T$ to $\nabla u$ combining the weak convergence with the convergence
of the norm.
Since $u_T$ is constant on each inclusion $B_i$, by a trace estimate on $B_i$ we have
$$
\expec{ u_T^2 \mathds 1_{\calB}} \,\lesssim \, \expec{ u_T^2  \mathds1_{\R^d\setminus\calB} +|\nabla u_T|^2},
$$
so that also $\lim_{T\to \infty} \frac1T \expec{ u_T^2 \mathds 1_{\calB}}=0$, and therefore  $\lim_{T\to \infty} \frac1T \expec{ u_T^2 }=0$, as claimed.

\subsection{Proposition~\ref{prop:unif-uT}: Caflisch-Luke estimates}

We split the proof into two steps.
We first prove the bounds on $\nabla u_T$, and then turn to the bounds on $u_T$ itself.
All the multiplicative constants in this proof depend on $\rho$.

\medskip

\step1 Proof of \eqref{e.unif-uT-0}.

By Lemma~\ref{lem:massive}, $\sup_{x\in \R^d} \fint_{B(x)} \frac1T u_T^2+ |\nabla u_T|^2 \lesssim T^{1+\frac d2}$, so that
$ \fint_{B} u_T^2+ |\nabla u_T|^2 \in L^\infty(\Omega)$.
By \eqref{e.compact} in Lemma~\ref{lem:compactness}, for all $p\ge 1$, all $R\gg 1$, and all $\sqrt{T}\gg R$,
we have (by adding the last RHS term to both sides of the inequality)
\begin{multline*} 
\expec{\Big(\fint_{B^R} |\nabla u_T|^2\Big)^p}^\frac1p +\expec{\Big|\fint_{B^{2R}} \frac1{\sqrt T} u_T\Big|^{2p}}^\frac1p\\
\le \, C  \sum_{n=1}^N \expec{\Big|\int_{B^{2R}}  \nabla u_T\cdot g_{R,n}\Big|^{2p} }^\frac1p
+C  R^2 
+ C \expec{\Big|\fint_{B^{2R}} \frac1{\sqrt T} u_T\Big|^{2p}}^\frac1p.
\end{multline*}
By Proposition~\ref{prop:unif-nablauT}, this entails
\begin{multline*} 
\expec{\Big(\fint_{B^R} |\nabla u_T|^2\Big)^p}^\frac1p +\expec{\Big|\fint_{B^{2R}} \frac1{\sqrt T} u_T\Big|^{2p}}^\frac1p
 \\
\lesssim \, p^{\gamma} R^{2-d} \Big(1
+\frac1{T} \expec{\Big|\fint_{B^{2R}} \frac1{\sqrt T}u_T\Big|^{2p}}^\frac1p+  \expec{\Big(\fint_{B^{R}}  |\nabla u_T|^2\Big)^{p} }^\frac1p\Big)
+C  R^2 ,
\end{multline*}
so that one may absorb the first RHS sum into the LHS for $R \gg p^{\frac \gamma{d-2}}$ and $d>2$. 
This yields the desired estimate
$$
\expec{\Big(\fint_{B} |\nabla u_T|^2\Big)^p}^\frac1p \le R^d \expec{\Big(\fint_{B^R} |\nabla u_T|^2\Big)^p}^\frac1p \,\le\, Cp^{\gamma \frac{d+2}{d-2}} .
$$

\medskip

\step2 Proof of \eqref{e.unif-uT-1}.

For simplicity we assume in this step that Proposition~\ref{prop:unif-nablauT} holds for all $R\ge 1$ (in the general case, it is enough to 
replace $B$ below by $B^r$ for some $r\gg1$ sufficiently large).
By Poincar\'e's inequality on $B$, 
$$
\fint_{B} u_T^2 \,\lesssim\, \Big(\fint_{B} u_T \Big)^2+\fint_{B} |\nabla u_T|^2,
$$
so that by Step~1, stationarity of $u_T$, and the triangle inequality,
\begin{equation} \label{e.bd-uT-1}
\expec{\Big(\fint_B u_T^2\Big)^p}^\frac1p \, \lesssim \,\expec{\Big| \fint_{B^{\sqrt{T}}} u_T-\fint_{B}u_T \Big|^{2p}}^\frac1p+\expec{\Big| \fint_{B^{\sqrt{T}}} u_T \Big|^{2p}}^\frac1p
+p^{\gamma \frac{d+2}{d-2}} .
\end{equation}
We split the rest of the proof into three substeps: We first estimate the  second RHS of \eqref{e.bd-uT-1}, then the  first RHS term of \eqref{e.bd-uT-1}, 
and we finally conclude.

\medskip

\substep{2.1} Control of the second RHS term of \eqref{e.bd-uT-1}.

We appeal to Proposition~\ref{prop:unif-nablauT} for $R=\sqrt T/2$, which yields, in combination with Step~1,
\begin{equation*}
\frac1T \expec{\Big| \fint_{B^{\sqrt{T}}} u_T \Big|^{2p}}^\frac1p\,\lesssim\, p^\gamma \sqrt{T}^{2-d}\bigg( 1
+\frac1{T^2}  \expec{\Big| \fint_{B^{\sqrt{T}}} u_T \Big|^{2p}}^\frac1p
+p^{\gamma \frac{d+2}{d-2}} \bigg).
\end{equation*}
By the deterministic energy estimate \eqref{e.massive-apriori} in form of $\fint_{B^{\sqrt{T}}} u_T^2 \lesssim T^2$, we can control the second RHS term, and obtain
\begin{equation}\label{e.crux-base}
\expec{\Big| \fint_{B^{\sqrt{T}}} u_T \Big|^{2p}}^\frac1p\,\lesssim\, \sqrt{T}^{4-d}  p^{\gamma \frac {2d}{d-2}}  .
\end{equation}

\medskip

\substep{2.2} Control of the first RHS term of \eqref{e.bd-uT-1}.

For all $r\ge 1$, let $h_r$ denote the unique radial solution of $-\triangle h_r=\frac1{|B^r|} \mathds{1}_{B^r}-\frac1{|B|} \mathds{1}_B$.
Then %
$$
 \fint_{B^{\sqrt{T}}} u_T-\fint_{B}u_T \,=\, \int \nabla u_T \cdot g_{\sqrt T},
$$
with $g_{\sqrt T}=\nabla h_{\sqrt T}$. By solving the equation for $h_r$ in radial coordinates (see e.g. \cite[Proof of Theorem~2, Step~3]{GNO-quant}), we obtain that $\supp {g_{\sqrt{T}}} \subset B^{\sqrt T}$ and $|g_{\sqrt{T}}(x)|\lesssim (1+|x|)^{1-d}$.
We may therefore appeal to Corollary~\ref{cor:unif-nablauT} with $R=\sqrt T$
which yields in combination with Step~1
\begin{eqnarray}\label{e.crux}
\expec{\Big| \fint_{B^{\sqrt{T}}} u_T-\fint_{B}u_T \Big|^{2p}}^\frac1p &\lesssim& \mu_d(\sqrt T) \Big( p^{\gamma \frac {2d}{d-2}}
 +p^ \gamma  \frac1{T^2}  \expec{\Big| \fint_{B} u_T^2 \Big|^{p}}^\frac1{p}\Big).
\end{eqnarray}

\medskip

\substep{2.3} Proof of \eqref{e.unif-uT-1}.

\nopagebreak 
The combination of \eqref{e.bd-uT-1}, \eqref{e.crux-base},   the triangle inequality, and \eqref{e.crux}  yields
for all $p\ge 1$  
\begin{equation}\label{e.crux-again}
\expec{\Big(\fint_B u_T^2\Big)^p}^\frac1p \, \lesssim \,\mu_d(\sqrt T)  p^{\gamma \frac {2d}{d-2}}
  + p^\gamma \mu_d(\sqrt T)\frac1{T^2} \expec{\Big(\fint_B u_T^2\Big)^p}^\frac1p.
\end{equation}
It remains to absorb the last RHS term of \eqref{e.crux-again} into the LHS. To this aim we use the energy estimate in form of
$\fint_{B} u_T^2 \lesssim T^{2+\frac d2}$ in combination with the additional decaying factor $\frac1{T^2}$
to the effect that 
$$
 \mu_d(\sqrt T)\frac1{T^2} \expec{\Big(\fint_B u_T^2\Big)^p}^\frac1p \,\lesssim\, \expec{\Big(\fint_B u_T^2\Big)^{(1-\alpha_d) p}}^\frac1p ,
$$
with the notation
$$
\alpha_d:=\left\{ 
\begin{array}{rcl}
d=3&:&\frac{3}{7},\\
d=4&:&\frac12^-,\\
d>4&:&\frac{4}{d+4},
\end{array}
\right.
$$
and where $\frac12^-$ means any exponent strictly less than $\frac12$.

We then use Jensen's inequality and Young's inequality with exponents $(\frac1{\alpha_d},\frac{1}{1-\alpha_d})$ and get for all $C\gg 1$
$$
 p^\gamma \mu_d(\sqrt T)\frac1{T^2} \expec{\Big(\fint_B u_T^2\Big)^p}^\frac1p \,\lesssim_d\, C p^{\gamma \frac1{\alpha_d}} +\frac1C \expec{\Big(\fint_B u_T^2\Big)^p}^\frac1p.
$$
We may thus absorb the last RHS term into the LHS of \eqref{e.crux-again}.
The desired estimate   \eqref{e.unif-uT-1} follows.

%%%%%%%%%%%%%
%%%%%%%%%%%

\subsection{Proposition~\ref{prop:unif-nablauT}: Decay of averages of $(\frac1{\sqrt T}u_T,\nabla u_T)$}

Starting point is the $p$-version of the multiscale Poincar\'e inequality: by \cite[Proposition 1.10 (ii)]{DG2}, \eqref{e.SG} entails for all centered random variables $Y$ and all exponents $p\ge1$
\begin{equation}\label{e.SGp}
\expec{|Y|^{2p}}
\le (Cp^{2})^p\,\int_1^\infty\expec{\bigg(\int_{\R^d} \oscd{Y}{B^\ell(z)}dz\bigg)^p} \ell^{-dp}e^{-\frac1C \ell}d\ell,
\end{equation}
which we shall apply to the random variable $Y_R = \int (\frac1{\sqrt T}u_T,\nabla u_T) \cdot g_{R}$.
To control the RHS of \eqref{e.SGp}, we shall distinguish whether $1\ll \ell \le R$ or $\ell \ge R\gg 1$,
and for each regime we shall consider far-field and near-field contributions separately.
Let $\ell \ge 1$ and $z \in \R^d$.
Let $\calP'\in \calS_\rho$ be such that $\calP'|_{\R^d\setminus B^\ell(z)}=\calP|_{\R^d\setminus B^\ell(z)}$, and recall that $w_T=u_T'-u_T$
(where $u_T'$ is the solution associated with $\calP'$).
We then have
$$
\delta_{B^\ell(z)} Y_R \,=\, \int  (\frac1{\sqrt T} w_T,\nabla w_T)  \cdot g_R
\quad \text{and} \quad
\osc{Y_R}{B^\ell(z)} \,=\, \sup_{\calP'} \delta_{B^\ell(z)} Y_R.
$$
We  split the rest of the proof into 4 steps.
In Step~1, we establish a pointwise decay estimate for $w_T(x)$ far from the source term, that is, for $x \notin B^{2\ell}(z)$.
In Steps~2 and 3, we consider the regimes $\ell \ge R$ and $\ell \le R$, respectively. We conclude in Step~4.

\medskip

\step1 Preliminary estimate: For all $\ell \gg 1$, $z\in \R^d$, and $x\in \R^d\setminus B^{2\ell}(z)$, we have
\begin{multline}\label{e.morceau2}
\frac1{\sqrt T} |w_T(x)|+|\nabla w_T(x)|\,\lesssim\, (\ell+|x-z|)^{1-d}e^{-\frac{\ell+|x-z|}{C\sqrt T}}
\\
 \times\ell^\frac {(d-1)(d+2)}{2d}\Big( \ell^{d+2}+\ell^{2}\frac1{T^2} \int_{B^{2\ell}(z)} u_T^2+ \int_{B^{2\ell}(z)} |\nabla u_T|^2\Big)^\frac12.
\end{multline}
To prove this estimate, we first introduce the notation 
\begin{eqnarray*}
B_{\ell}^\calB(z)& := &B^\ell(z) \cup_{B\in \calB: B\cap B^\ell(z)\neq \emptyset}B \,\subset \,B^{\ell+\sqrt d}(z),
\\
\bar B_\ell^\calB(z)&:=&\{x\in \R^d\,|\, \dist(x,B_{\ell}^\calB(z))<1\}\,\subset\, B^{\ell+\sqrt d +1}(z),
\\
\calB_{\ell,z}&:=&\{B\in \calB\,|\, B\cap B^\ell(z)= \emptyset\}.
\end{eqnarray*}
Set $w_T^{\ell,z}:=w_T \chi_{\ell,z}$, where  $\chi_{\ell,z}:\R^d \to [0,1]$ is a smooth cut-off for $B_{\ell}^\calB(z)$ in $\bar B_\ell^\calB(z)$
(that is, such that $\chi_{\ell,z}|_{B_\ell^\calB(z)}\equiv 1$ and $\chi_{\ell,z}|_{\R^d \setminus \bar B_\ell^\calB(z)}\equiv 0$).
By definition,  $w_T$ satisfies 
\begin{multline} \label{a.10}
\frac1T (w_T-w_T^{\ell,z})-\triangle (w_T-w_T^{\ell,z})\,=\, \nabla w_T \cdot \nabla \chi_{\ell,z}+ \nabla \cdot (w_T \nabla  \chi_{\ell,z})
\quad \text{ in } \R^d\setminus (\cup_{B\in \calB_{\ell,z}}B), 
\\
\forall i \text{ such that }B_i \in \calB_{\ell,z}:\quad w_T|_{B_i} \equiv w_{T,i} \in \R,  \quad \int_{\partial B_i} \partial_n w_T = 0.
\end{multline}
Denote by $G_T$ the Green's function of Lemma~\ref{lem:Green} associated with the point set $\calP_{\ell,z}:=\{x \,|\, x \in \calP,B_x \cap B^\ell(z) =\emptyset\}$.
By the choice of the cut-off, the Green representation formula yields
for all $x\in \R^d\setminus B^{\ell+\sqrt d+1}(z)$  
\begin{eqnarray*}
\nabla w_T(x)&=&\nabla_x(w_T(x)-w_T^{\ell,z}(x))\\
&=&2 \int_{\R^d} \nabla_x G_T(x,y)(\nabla w_T \cdot \nabla \chi_{\ell,z})(y)dy
+\int_{\R^d} \nabla_x G_T(x,y) \triangle \chi_{\ell,z}(y) w_T(y)dy,
\\
w_T(x)&=&2\int_{\R^d}  G_T(x,y)(\nabla w_T \cdot \nabla \chi_{\ell,z})(y)dy+\int_{\R^d} G_T(x,y) \triangle \chi_{\ell,z}(y) w_T(y)dy,
\end{eqnarray*}
so that by \eqref{e.Green2} \&~\eqref{e.Green1} followed by the Gagliardo-Nirenberg-Sobolev inequality ($d>2$) and the Cauchy-Schwarz inequality,
\begin{eqnarray*}
\lefteqn{\frac1{\sqrt{T}} |w_T(x)|+|\nabla_x w_T(x)|}
\\
&\lesssim & (1+\dist(x,B^{\ell+\sqrt d+1}(z))^{1-d}e^{-\frac{\dist(x,B^{\ell+\sqrt d+1}(z))}{C\sqrt T}} 
\int_{B^{\ell+\sqrt d+1}(x)\setminus B^{\ell-\sqrt d-1}(x)} |\nabla w_T|+|w_T|
\\
&\lesssim& (1+\dist(x,B^{\ell+\sqrt d+1}(z))^{1-d}e^{-\frac{\dist(x,B^{\ell+\sqrt d+1}(z))}{C\sqrt T}}
\ell^\frac {(d-1)(d+2)}{2d} \Big(\int_{\R^d} |\nabla w_T|^2\Big)^\frac12.
\end{eqnarray*}
Combined with Corollary~\ref{cor:apriori}, this finally implies the desired decay estimate of $w_T$. 

\medskip

\step2 Estimate of $\osc{Y_R}{B^\ell(z)}$ for $\ell\ge R\gg 1$.

We start with the  near-field contribution, that is, for $|z|\lesssim \ell$.
In this case we have by Cauchy-Schwarz' inequality and the support condition on $g$
\begin{equation}\label{e.key-1}
\Big(\int  (\frac1{\sqrt T}w_T,\nabla w_T)  \cdot g_R\Big)^2\,\lesssim \, R^{-d} \int_{\R^d} \frac1T w_T^2+ |\nabla w_T|^2.
\end{equation}
Since $\ell \ge \rho \ge \rho_\m$, by \eqref{e.diff-uT2} 
\begin{equation}\label{e.key-1.1}
\int_{\R^d} \frac1T w_T^2+|\nabla w_T|^2 \,\le\,  C \ell^{d+2}+C  \ell^{2} \frac1{T^2} \int_{B^{2\ell}(z)} u_T^2
+ \int_{B^{2\ell}(z)} |\nabla u_T|^2 .
\end{equation}
Hence, \eqref{e.key-1} turns into
\begin{equation}\label{e.key-2}
\oscd{Y_R}{B^\ell(z)}\,\le \, C R^{-d}  \Big( \ell^{d+2}+\ell^{2} \frac1{T^2} \int_{B^{2\ell}(z)} u_T^2
+ \int_{B^{2\ell}(z)} |\nabla u_T|^2\Big).
\end{equation}
We then turn to the far-field contribution, that is, for $|z| \gg  \ell$.
In this case we apply \eqref{e.morceau2}  for $x\in B^{2R}$ and $|z|\gg \ell \ge R$ (so that $x\notin B^{2\ell}(z)$)
in form of 
\begin{equation}\label{e.key-3}
\oscd{Y_R}{B^\ell(z)} \,\lesssim \,  \ell^\frac{(d-1)(d+2)}{d} \Big( \ell^{d+2}+\ell^{2} \frac1{T^2} \int_{B^{2\ell}(z)} u_T^2+ \int_{B^{2\ell}(z)} |\nabla u_T|^2\Big) \frac{e^{-\frac{|z|}{C\sqrt T}}}{|z|^{2(d-1)}}.
\end{equation}

\medskip

\step3 Estimate of $\osc{Y_R}{B^\ell(z)}$ for $R \ge \ell\gg 1$.

We first use that $|g_R|\lesssim R^{-d} \mathds 1_{B^{2R}}$ to the effect of 
$$
\Big(\int  (\frac1{\sqrt T}w_T,\nabla w_T)  \cdot g_R\Big)^2\,\lesssim \, R^{-2d} \Big(\int_{B^{2R}} \frac{1}{\sqrt T}|w_T|+|\nabla w_T| \Big)^2.
$$
We start with the far-field contribution $|z|\gg R$ for which \eqref{e.morceau2} yields after integration over $x \in B^{2R}$
\begin{equation}\label{e.key-4}
\oscd{Y_R}{B^\ell(z)} \,\lesssim\,  R^{-2(d-1)}\ell^\frac {(d-1)(d+2)}{d} \Big( \ell^{d+2}+\ell^{2} \frac1{T^2} \int_{B^{2\ell}(z)} u_T^2+ \int_{B^{2\ell}(z)} |\nabla u_T|^2\Big) \frac{e^{-\frac{|z|}{C\sqrt T}}}{|z|^{2(d-1)}}.
\end{equation}
We then turn to the near-field contribution $|z|\lesssim R$.
Let $n \in \N$ be the smallest integer such that $B^{2R} \subset B^{2^{n+1}\ell}(z)$, and note that $n\lesssim \log (\frac {R}{\ell} +2)$.
We bound the integral on $B^{2R}$ of non-negative integrands as 
$$
\int_{B^{2R}}\,\le\, \int_{B^{2\ell}(z)}+\sum_{i=1}^{n} \int_{B^{2^{i+1}\ell}(z)\setminus B^{2^{i}\ell}(z)}.
$$
On the first set, we use Cauchy-Schwarz' inequality together with \eqref{e.key-1.1}
to the effect of 
\begin{equation}\label{e.key-5.1}
\int_{B^{2\ell}(z)} \frac{1}{\sqrt T}|w_T|+ |\nabla w_T| \,\lesssim\, \ell^\frac{d}2 \Big( \ell^{d+2}+\ell^{2} \frac1{T^2} \int_{B^{2\ell}(z)} u_T^2
+ \int_{B^{2\ell}(z)} |\nabla u_T|^2 \Big)^\frac12.
\end{equation}
Since the other sets satisfy $x\in B^{2^{i+1}\ell}(z)\setminus B^{2^{i}\ell}(z) \,\implies \, x \notin B^{2\ell}(z)$, 
one may appeal to \eqref{e.morceau2}, which yields 
\begin{multline}\label{e.key-5.2}
\int_{B^{2^{i+1}\ell}(z)\setminus B^{2^{i}\ell}(z) } \frac{1}{\sqrt T}|w_T|+ |\nabla w_T| \,\lesssim\, (2^i\ell)^d (2^i\ell)^{1-d}  
\\
\times \ell^\frac {(d-1)(d+2)}{2d}\Big( \ell^{d+2}+\ell^{2}\frac1{T^2} \int_{B^{2\ell}(z)} u_T^2+ \int_{B^{2\ell}(z)} |\nabla u_T|^2\Big)^\frac12.
\end{multline}
Summing \eqref{e.key-5.1} and \eqref{e.key-5.2} over $i=1,\dots, n$ then entails 
\begin{equation*} 
\int_{B^{2R}} \frac{1}{\sqrt T}|w_T|+ |\nabla w_T| \,\lesssim\, 
 R \ell^\frac {(d-1)(d+2)}{2d}\Big( \ell^{d+2}+\ell^{2}\frac1{T^2} \int_{B^{2\ell}(z)} u_T^2+ \int_{B^{2\ell}(z)} |\nabla u_T|^2\Big)^\frac12,
\end{equation*}
so that for $|z| \lesssim R$,
\begin{equation}\label{e.key-5}
\oscd{Y_R}{B^\ell(z)} \,\lesssim\,  R^{-2(d-1)}\ell^\frac {(d-1)(d+2)}{d} \Big( \ell^{d+2}+\ell^{2} \frac1{T^2} \int_{B^{2\ell}(z)} u_T^2+ \int_{B^{2\ell}(z)} |\nabla u_T|^2\Big) .
\end{equation}

\step4 Proof of \eqref{e.unif-nablauT}.

\nopagebreak

Starting point is \eqref{e.SGp}, which we split into the two contributions $\ell \le R$ and $\ell\ge R$:
\begin{multline*} 
\expec{|Y_R|^{2p}} \,
\le \;(Cp^{2})^p\,\int_1^R\expec{\bigg(\int_{\R^d} \oscd{Y_R}{B^\ell(z)}dz\bigg)^p} \ell^{-dp}e^{-\frac1C \ell}d\ell
\\
+(Cp^{2})^p\,\int_R^\infty\expec{\bigg(\int_{\R^d} \oscd{Y_R}{B^\ell(z)}dz\bigg)^p} \ell^{-dp}e^{-\frac1C \ell}d\ell.
\end{multline*}
Since $\ell \mapsto  \oscd{Y_R}{B^\ell(z)}$ is non-decreasing, one may assume wlog that $\ell\gg 1$, in which 
case the estimates of Steps~2 and~3 are in force.
For $\ell \le R$, we average \eqref{e.key-4} and \eqref{e.key-5} on balls of size $R$ to the effect that 
\begin{eqnarray*}
\lefteqn{\int_{\R^d} \oscd{Y_R}{B^\ell(z)}dz}
\\
 &=& \int_{\R^d} \fint_{B^R(z)} \oscd{Y_R}{B^\ell(z')}dz'dz
\\
&\lesssim&
\int_{B^R} R^{-2(d-1)}\ell^\frac {(d-1)(d+2)}{d} \ell^d \Big( \ell^{2}+\ell^{2} \frac1{T^2} \fint_{B^{2R}} u_T^2+ \fint_{B^{2R}} |\nabla u_T|^2\Big)
\\
&&+\int_{\R^d \setminus B^R} R^{-2(d-1)}\ell^\frac {(d-1)(d+2)}{d} \ell^d \Big( \ell^{2}+\ell^{2} \frac1{T^2} \fint_{B^{2R}(z)} u_T^2+  \fint_{B^{2R}(z)} |\nabla u_T|^2\Big) \frac{e^{-\frac{|z|}{C\sqrt T}}}{|z|^{2(d-1)}}dz.
\end{eqnarray*}
Combined with the triangle inequality for $\int_1^R \expec{|\cdot|}$, this yields for  $\alpha=d+2+\frac {(d-1)(d+2)}{d} $ by stationarity of $u_T$ and $\nabla u_T$
and using that $2(d-1)>d$ for $d>2$ to treat the integral over $\R^d \setminus B^R$, and Jensen's inequality to pass from $\fint_{B^{2R}}$ to $\fint_{B^R}$,
\begin{eqnarray}
\lefteqn{\Big(\int_1^R\expec{\bigg(\int_{\R^d} \oscd{Y_R}{B^\ell(z)}dz\bigg)^p} \ell^{-dp}e^{-\frac1C \ell}d\ell\Big)^\frac1p}\nonumber
\\
&\lesssim& R^{2-d} \expec{1+\Big(\fint_{B^{2R}} \frac1{T^2}u_T^2\Big)^{p}+\Big(\fint_{B^{2R}} |\nabla u_T|^2\Big)^p}^\frac1p 
\int_1^R  \ell^\alpha e^{-\frac\ell {Cp}} d\ell \qquad\qquad\qquad\qquad\nonumber
\\
&\lesssim& R^{2-d} p^{\alpha+1}\Big(1+\frac1{T}\expec{\Big(\fint_{B^{2R}} \frac1{\sqrt T} u_T\Big)^p}^\frac1p+(1+\frac{R^2}{T^2})\expec{\Big(\fint_{B^{R}} |\nabla u_T|^2\Big)^p}^\frac1p\Big),
\label{e.key-7}
\end{eqnarray}
where we used Poincar\'e's inequality on $B^R$ in the last line.
For $\ell \ge R$, we integrate  \eqref{e.key-2} on $B^\ell$ and \eqref{e.key-3} on $\R^d \setminus B^\ell$, which yields
\begin{eqnarray*}
\lefteqn{\int_{\R^d} \oscd{Y_R}{B^\ell(z)}dz}
\\
&\lesssim&
C R^{-d} \ell^{2d+2}  \Big( 1+ \frac1{T^2} \fint_{B^{3\ell}(0)} u_T^2
+ \fint_{B^{3\ell}(0)} |\nabla u_T|^2\Big)
\\
&&+\int_{\R^d \setminus B^R} \ell^{d+2+\frac {(d-1)(d+2)}{d}} \Big( 1+ \frac1{T^2} \fint_{B^{2\ell}(z)} u_T^2+ \fint_{B^{2\ell}(z)} |\nabla u_T|^2\Big) \frac{e^{-\frac{|z|}{C\sqrt T}}}{|z|^{2(d-1)}}dz.
\end{eqnarray*}
As above, combined with the triangle inequality for $\int_R^\infty \expec{|\cdot|}$, this yields for  $\alpha=d+2+\frac {(d-1)(d+2)}{d} $ by stationarity of $u_T$ and $\nabla u_T$
and using that $2(d-1)>d$ for $d>2$ to treat the integral over $\R^d \setminus B^R$, and Jensen's inequality to pass from $\fint_{B^{2\ell}}$ to $\fint_{B^{2R}}$ (and $\fint_{B^{R}}$)
since $\ell \ge R$,
\begin{eqnarray}
\lefteqn{\Big(\int_R^\infty\expec{\bigg(\int_{\R^d} \oscd{Y_R}{B^\ell(z)}dz\bigg)^p} \ell^{-dp}e^{-\frac1C \ell}d\ell\Big)^\frac1p}\nonumber
\\
&\lesssim& R^{2-d} \expec{1+\Big(\frac1{T^2}\fint_{B^{2R}}u_T^2\Big)^p+\Big(\fint_{B^{2R}} |\nabla u_T|^2\Big)^p}^\frac1p 
\int_1^\infty  \ell^\alpha e^{-\frac\ell {Cp}} d\ell \qquad\qquad\qquad\qquad \nonumber
\\
&\lesssim& R^{2-d} p^{\alpha+1}\Big(1+\frac1{T}\expec{\Big(\fint_{B^{2R}}\frac1{\sqrt{T}}u_T\Big)^{2p}}^\frac1p+(1+\frac{R^2}{T^2})\expec{\Big(\fint_{B^{R}} |\nabla u_T|^2\Big)^p}^\frac1p\Big),
\label{e.key-8}
\end{eqnarray}
where we also used Poincar\'e's inequality on $B^{2R}$ in the last line.
The combination of \eqref{e.key-7} and  \eqref{e.key-8} with \eqref{e.SGp} concludes the proof.

\subsection{Proof of Corollary~\ref{cor:unif-nablauT}: Finer decay of averages}

The proof of this corollary has the same structure as the proof of Proposition~\ref{prop:unif-nablauT}.
From a technical point of view we cannot take advantage any longer of the scaling of $g_R$ wrt $R$ (cf. $g_R=R^{-d}g(\frac \cdot R)$
versus $g_R(x)=(1+|x|)^{1-d} \mathds1_{B^R}$).
In terms of estimates, there is no difference when we use the Cauchy-Schwarz inequality (as in Step~2) to control
terms of the form $\int \nabla w_T \cdot g_R$.
The only significant difference is that $L^\infty$ bounds on $g_R$ now yield suboptimal scalings, which
compels us to be more precise to unravel cancellations.  This concerns Step~3, which
we therefore presently adapt in detail, before concluding as before.

\medskip

Set $Y_R=\int \nabla u_T \cdot g_R$, and consider the regime  $R \ge \ell\gg 1$.
The aim is to estimate $\osc{Y_R}{B^\ell(z)}$.
We start with the  far-field contribution $|z|\gg R$, which is easier.
In that case, by definition of $g_R$,
\begin{eqnarray*}
\int |\nabla w_T| |g_R|&\le & \sup_{B^R}|\nabla w_T| \int_{B^R} (1+|x|)^{1-d} dx 
\\
&\lesssim & R  \sup_{B^R}|\nabla w_T|,
\end{eqnarray*}
so that  by \eqref{e.morceau2}  we obtain
\begin{equation}\label{cri:key-4}
\oscd{Y_R}{B^\ell(z)} \,\lesssim\,   R^2 \ell^\frac {(d-1)(d+2)}{d} \Big( \ell^{d+2}+\ell^{2} \frac1{T^2} \int_{B^{2\ell}(z)} u_T^2+ \int_{B^{2\ell}(z)} |\nabla u_T|^2\Big) \frac{e^{-\frac{|z|}{C\sqrt T}}}{|z|^{2(d-1)}}.
\end{equation}
We then turn to the near-field contribution $|z|\lesssim R$.
Let $n \in \N$ be the smallest integer such that $B^{R} \subset B^{2^{n+1}\ell}(z)$, and note that $n\lesssim \log (\frac {R}{\ell} +2)$.
Define $C^0=B^{2\ell}(z)$, $C^i=B^{2^{i+1}\ell}(z)\setminus B^{2^{i}\ell}(z)$ for all $1\le i\le n$.
We bound the integral on $B^{R}$ of non-negative integrands as 
$$
\int_{B^{R}}\,\le\, \int_{C^0}+\sum_{i=1}^n \int_{C^i}.
$$
On the one hand,  by  \eqref{e.morceau2}, we have for  $1\le i\le n$,
\begin{equation*}
\sup_{C^i}  |\nabla w_T| \,\lesssim \, (2^i\ell)^{1-d}  
 \ell^\frac {(d-1)(d+2)}{2d}\Big( \ell^{d+2}+\ell^{2}\frac1{T^2} \int_{B^{2\ell}(z)} u_T^2+ \int_{B^{2\ell}(z)} |\nabla u_T|^2\Big)^\frac12,
 \end{equation*}
so that
\begin{eqnarray}
\lefteqn{\sum_{i=1}^n \int_{C^i} |\nabla w_T||g_R|}\nonumber
\\ &\lesssim &  \ell^\frac {(d-1)(d+2)}{2d}\Big( \ell^{d+2}+\ell^{2}\frac1{T^2} \int_{B^{2\ell}(z)} u_T^2+ \int_{B^{2\ell}(z)} |\nabla u_T|^2\Big)^\frac12\nonumber
\\
&&\times \sum_{i=1}^n \int_{C^i} |x-z|^{1-d} (1+|x|)^{1-d}dx\nonumber
\\
&\lesssim &(1+|z|)^{2-d} \ell^\frac {(d-1)(d+2)}{2d}\Big( \ell^{d+2}+\ell^{2}\frac1{T^2} \int_{B^{2\ell}(z)} u_T^2+ \int_{B^{2\ell}(z)} |\nabla u_T|^2\Big)^\frac12.\label{cri:001}
\end{eqnarray}
On the other hand, by  \eqref{e.key-1.1}, for $i=0$ we have
$$
\int_{C^0} |\nabla w_T|^2 \,\lesssim \,   \ell^{d+2}+\ell^{2}\frac1{T^2} \int_{B^{2\ell}(z)} u_T^2+ \int_{B^{2\ell}(z)} |\nabla u_T|^2.
$$
If  $|z|\le 4\ell$ we thus obtain by Cauchy-Schwarz' inequality on $C^0$ 
\begin{eqnarray}
\lefteqn{\int_{C^0} |\nabla w_T||g_R|}\nonumber 
\\ &\lesssim &    \Big( \ell^{d+2}+\ell^{2}\frac1{T^2} \int_{B^{2\ell}(z)} u_T^2+ \int_{B^{2\ell}(z)} |\nabla u_T|^2\Big)^\frac12,\nonumber
\\
&&\times \Big(\int_{B^{2\ell}(z)}  (1+|x|)^{2(1-d)}dx\Big)^\frac12\nonumber
\\
&\lesssim & \mathds{1}_{|z|\le 4\ell}    \Big( \ell^{d+2}+\ell^{2}\frac1{T^2} \int_{B^{2\ell}(z)} u_T^2+ \int_{B^{2\ell}(z)} |\nabla u_T|^2\Big)^\frac12,\label{cri:002}
\end{eqnarray}
whereas of $|z|>4\ell$, $\sup_{C^0} |g_R|\lesssim (1+|z|)^{1-d}$ and by Cauchy-Schwarz' inequality on $C^0$  again
\begin{eqnarray}
\int_{C^0} |\nabla w_T||g_R| &\lesssim & (1+|z|)^{1-d} \int_{B^{2\ell}(z)} |\nabla w_T|
\nonumber
\\
&\lesssim & (1+|z|)^{1-d}  \ell^\frac{d}2    \Big( \ell^{d+2}+\ell^{2}\frac1{T^2} \int_{B^{2\ell}(z)} u_T^2+ \int_{B^{2\ell}(z)} |\nabla u_T|^2\Big)^\frac12. \label{cri:003}
\end{eqnarray}
Estimates~\eqref{cri:001}--\eqref{cri:003}  combine to 
\begin{equation*}
\int |\nabla w_T| |g_R|\,\lesssim\, (\mathds{1}_{|z|\le 4\ell}+(1+|z|)^{2-d} \ell^\frac {(d-1)(d+2)}{2d})   \Big( \ell^{d+2}+\ell^{2}\frac1{T^2} \int_{B^{2\ell}(z)} u_T^2+ \int_{B^{2\ell}(z)} |\nabla u_T|^2\Big)^\frac12,
\end{equation*}
from which we infer that we have in the regime  
$|z|\lesssim R$ 
\begin{multline}\label{cri:key-5}
\oscd{Y_R}{B^\ell(z)} \,\lesssim\,   (\mathds{1}_{|z|\le 4\ell}+( |z|+1)^{2(2-d)} \ell^{ \frac {(d-1)(d+2)}{d}})
\\
\times \Big( \ell^{d+2}+\ell^{2}\frac1{T^2} \int_{B^{2\ell}(z)} u_T^2+ \int_{B^{2\ell}(z)} |\nabla u_T|^2\Big).
\end{multline}

\medskip

We conclude with the adaptation of Step~4, and only treat the contribution $1\le \ell \le R$, which is the dominating one. 
As opposed to Step~4, we \emph{do not} average \eqref{cri:key-4} and \eqref{cri:key-5} on balls of size $R$, 
and obtain
\begin{eqnarray*}
\lefteqn{\int_{\R^d} \oscd{Y_R}{B^\ell(z)}dz}
\\
 &\lesssim& \ell^{d+2}\int_{B^{4\ell}} \Big( 1+\frac1{T^2} \fint_{B^{2\ell}(z)} u_T^2+ \fint_{B^{2\ell}(z)} |\nabla u_T|^2\Big)dz
\\
&&+ \ell^{\frac {(d-1)(d+2)}{d}}\ell^{d+2} \int_{B^R}( |z|+1)^{2(2-d)} \Big( 1+\frac1{T^2} \fint_{B^{2\ell}(z)} u_T^2+ \fint_{B^{2\ell}(z)} |\nabla u_T|^2\Big)dz
\\
&&+R^{2}\ell^\frac {(d-1)(d+2)}{d}\ell^{d+2} \int_{\R^d \setminus B^R} \Big( 1+ \frac1{T^2} \fint_{B^{2\ell}(z)} u_T^2+ \fint_{B^{2\ell}(z)} |\nabla u_T|^2\Big) \frac{e^{-\frac{|z|}{C\sqrt T}}}{|z|^{2(d-1)}}dz.
\end{eqnarray*}
Combined with the triangle inequality for $\int_1^R \expec{|\cdot|}$, this yields for  $\alpha=d+2+\frac {(d-1)(d+2)}{d} $ by stationarity of $u_T$ and $\nabla u_T$, and Jensen's inequality
\begin{eqnarray*}
\lefteqn{\Big(\int_1^R\expec{\bigg(\int_{\R^d} \oscd{Y_R}{B^\ell(z)}dz\bigg)^p} \ell^{-dp}e^{-\frac1C \ell}d\ell\Big)^\frac1p}\nonumber
\\
 &\lesssim& \Big(\int_{B^R}( |z|+1)^{2(2-d)}dz +R^2   \int_{\R^d \setminus B^R}  \frac{e^{-\frac{|z|}{C\sqrt T}}}{|z|^{2(d-1)}}dz\Big)  
 \\
 &&\times \expec{1+\Big(\frac1{T^2}\fint_{B^{2}}u_T^2\Big)^p+\Big(\fint_{B^{2}} |\nabla u_T|^2\Big)^p}^\frac1p 
 \int_1^R  \ell^\alpha e^{-\frac\ell {Cp}} d\ell
 \\
 &\lesssim& \mu_d(R) p^{\alpha+1}\expec{1+\Big(\frac1{T^2}\fint_{B^{2}}u_T^2\Big)^p+\Big(\fint_{B^{2}} |\nabla u_T|^2\Big)^p}^\frac1p .
\end{eqnarray*}
The desired estimate \eqref{e.unif-nablauT-crit} then follows.

\section{Proofs of the regularity results}\label{sec:regularity-results}

\subsection{Lemma~\ref{lem:reg1}: Quantitative approximate radiality}

We split the proof into three steps.
In the first step we lift the function $v$ using the solution $v^\infty_T$ involving an explicit radial solution on the whole space.  
In the second step, we prove that the difference $v-v^\infty_T$ satisfies a hole-filling estimate at the origin of the decaying factor  $\rho^{-\alpha}$. We then conclude 
in the third step by controlling the $L^2$-norm of $\nabla v^\infty_T$ using an energy estimate and the explicit radial solution.

\medskip

\step1 Lifting.

First observe that up to replacing $g_2$ by $g_2-\frac{\bar v}{T}$, one can assume that $\bar v=0$, which we do in the rest of the proof.
Consider the radial solution $v^\infty:r \mapsto -g_2 \frac{r^2}{2d}+(-g_1+\frac{g_2}d)\frac{2-d}{r^{d-2}}+ \frac{g_2}{2d}+(d-2)(-g_1+\frac{g_2}d)$ of
the equation $-\triangle v^\infty\,=\, g_2 \text{ in }\R^d\setminus B^1, \quad v^\infty\equiv 0 \text{ on }\partial B^1, \quad \fint_{\partial B^1} \nabla v^\infty\cdot n=g_1$
on $\R^d$.
We define $v^\infty_T$ as the unique solution in $H^1(B^\rho)$ of  
$\frac1T v^\infty_T-\triangle v^\infty_T\,=\, g_2 \text{ in }B^\rho \setminus B^1, \quad v_T^\infty\equiv 0 \text{ on }\partial B^1, \quad \fint_{\partial B^1} \nabla v_T^\infty\cdot n=g_1$, and
$v^\infty_T|_{\partial B^\rho} \equiv v^\infty|_{\partial B^\rho}$.
By uniqueness, $v^\infty_T$ is radial, so that $\int_{\partial B^1} |\nabla v^\infty_T\cdot n-g_1|=0$.
We then set $w:=v-v^\infty_T$, and note that $w$ satisfies
$$
\frac1T w-\triangle w\,=\, 0 \text{ in }B^\rho\setminus B^1, \quad w\equiv 0 \text{ on }\partial B^1, \quad \fint_{\partial B^1} \nabla w\cdot n=0
$$
and that we have
$$
\int_{\partial B^1} |\nabla v\cdot n-g_1|\,\le\,\int_{\partial B^1} |\nabla v_T^\infty\cdot n-g_1|+\int_{\partial B^1} |\nabla w\cdot n| \,=\, \int_{\partial B^1} |\nabla w\cdot n|.
$$
By elliptic regularity up to the boundary and trace estimates, this yields
\begin{equation}\label{e.reg-1:elliptic}
\int_{\partial B^1} |\nabla v\cdot n-g_1|\,\lesssim\, \Big(\int_{B^2\setminus B^1} |\nabla w|^2\Big)^\frac12 = \Big(\int_{B^2} |\nabla w|^2\Big)^\frac12
\end{equation}
after extending $w$ by 0 in $B^1$.

\medskip

\step2 Hole-filling argument.

Let $r\ge 2$, let $N \in \N$ to be fixed later, and let $\eta_r$ be a cut-off for $B^r$ in $B^{Nr}$ such that $|\nabla \eta_r|\lesssim \frac1{Nr}$ (in particular, $\eta_r\equiv 1$ on $B^1\subset B^{r}$).
We then proceed to the Caccioppoli argument and test the defining equation for $w$ by $\eta_r^2( w- c)$ for  $c = (\int_{B^{Nr} \setminus B^1} \eta_r^2)^{-1}\int_{B^{Nr} \setminus B^1} \eta_r^2 w$.
This yields after standard calculations (there is no boundary contribution on $\partial B^1$ since $\eta_\rho^2( w- c)$ is constant on $\partial B^1$ and $\int_{\partial B^1} \nabla w \cdot n=0$)
$$
\frac1T \int_{\R^d\setminus B^1} \eta_r^2 (w-c)^2 + \int_{\R^d \setminus B^1}  |\nabla (\eta_r (w-c))|^2 = \int_{\R^d \setminus B^1} |\nabla \eta_r|^2 (w-c)^2 .
$$
By the choice of the cut-off, this yields for all $\tilde c \in \R$,
$$
\int_{B^{r}}  |\nabla w|^2 \le \frac C{N^2r^2}  \int_{B^{Nr}\setminus B^{r}}  (w-c)^2  \,\le \,  \frac C{N^2r^2} \int_{B^{Nr}\setminus B^{r}}  (w-\tilde c)^2+\frac C{N^2r^2} (Nr)^d (\tilde c-c)^2,
$$
where $C$ might change from line to line
but only depends on $d$. 
It remains to choose $\tilde c$. Let $N\gg 1$ (and therefore $\rho\gg 1$) be so large that there exists $\kappa \sim_d 1$ (bounded uniformly in $N$) so that $\int_{B^{\kappa r} \setminus B^r}\eta_r^2=\int_{B^r \setminus B^1}\eta_r^2$.
We then set 
$$
\tilde c \,:=\, \frac{\int_{B^{Nr} \setminus B^r} \eta_r^2 w+\int_{B^{\kappa r} \setminus B^r} \eta_r^2 w}{\int_{B^{Nr} \setminus B^1} \eta_r^2}=\frac{\int_{B^{Nr} \setminus B^r} \eta_r^2 w+\int_{B^{\kappa r} \setminus B^r} \eta_r^2 w}{\int_{B^{Nr} \setminus B^r} \eta_r^2+\int_{B^{\kappa r} \setminus B^r} \eta_r^2}.
$$
On the one hand, by Poincar\'e's inequality on $B^{Nr} \setminus B^r$, 
$$
\frac C {N^2r^2} \int_{B^{Nr}\setminus B^{r}}  (w-\tilde c)^2 \, \le \, C\int_{B^{Nr}\setminus B^{r}}  |\nabla w|^2.
$$
On the other hand, by Poincar\'e's inequality on $B^{\kappa r}$,
\begin{eqnarray*}
(c-\tilde c)^2 &=& (\frac{\int_{B^{\kappa r} \setminus B^r} \eta_r^2 w-\int_{B^{r} \setminus B^1} \eta_r^2 w}{\int_{B^{Nr} \setminus B^1} \eta_r^2})^2
\\
&=&  (\frac{\int_{B^{\kappa r} \setminus B^r} \eta_r^2 (w-\frac{1}{\int_{B^{r} \setminus B^1} \eta_r^2}\int_{B^{r} \setminus B^1} \eta_r^2 w)}{\int_{B^{Nr} \setminus B^1} \eta_r^2})^2
\\
&\lesssim& \frac{1}{(Nr)^{2d}} (\kappa r)^d \int_{B^{\kappa r}} (w-\frac{1}{\int_{B^{r} \setminus B^1} \eta_r^2}\int_{B^{r} \setminus B^1} \eta_r^2 w)^2
\\
&\lesssim&  \frac{1}{(Nr)^{2d}} (\kappa r)^{d+2}  \int_{B^{\kappa r}} |\nabla w|^2 \\
&\stackrel{\kappa \sim 1}\lesssim & \frac1{N^{2+d}} \frac{1}{(Nr)^d} N^2r^2  \int_{B^{N r}} |\nabla w|^2.
\end{eqnarray*}
Gathering these estimates then yields
$$
\int_{B^{r}}  |\nabla w|^2\, \le\, C\int_{B^{Nr}\setminus B^{r}}  |\nabla w|^2+ \frac{C}{N^{2+d}} \int_{B^{N r}} |\nabla w|^2.
$$
Adding $C$ times the LHS to this inequality entails 
$$
\int_{B^{r}}  |\nabla w|^2 \,\le\, \frac{C}{C+1}(1+\frac1{N^{2+d}}) \int_{B^\rho }  |\nabla w|^2 \, \le \,\frac{C}{C+1/2}  \int_{B^{Nr} }  |\nabla w|^2
$$
provided $N \ge (2C+1)^\frac{1}{2+d}$.
This estimate can be iterated, which yields the hole-filling estimate for $\theta = \log(\frac{C+1/2}C)/\log N>0$ and $\rho \gg 1$
\begin{equation}\label{e.reg-1:hole}
\int_{B^2} |\nabla w|^2 \,\lesssim \, \rho^{-\theta} \int_{B^\rho} |\nabla w|^2
\end{equation}
(the multiplicative constant depends on $N$).

\medskip

\step3 Conclusion.

\nopagebreak
The combination of \eqref{e.reg-1:elliptic} and \eqref{e.reg-1:hole} yields for all $\rho \ge 4$
$$
\int_{\partial B^1} |\nabla v\cdot n-g_1| \le C \rho^{-\theta/2}  \Big(\int_{B^\rho} |\nabla w|^2\Big)^\frac12.
$$
The claim follows from the triangle inequality applied to $w=v-v^\infty_T$ provided we prove the energy estimate 
\begin{equation}\label{e.reg-1:energy}
\int_{B^\rho} |\nabla v^\infty_T|^2 \,\lesssim\, \rho^{d}(1+\frac{\rho^2}{T})(g_1^2+\rho^2 g_2^2) . 
\end{equation}
This is a direct consequence of the explicit formula for $v^\infty$ and of the following energy estimate for $v_T^\infty-v_\infty$:
$$
\frac1T \int_{B^\rho} (v_T^\infty-v^\infty)^2+ \int_{B^\rho} |\nabla (v_T^\infty-v_\infty)|^2 \,\lesssim\,\frac1T \int_{B^\rho} (v^\infty)^2.
$$

\subsection{Lemma~\ref{lem:Green}: Green's functions estimates}

We split the proof into five steps, and drop the subscripts $T$ for readability. 
Since the estimates are uniform with respect to the point sets in the class $\calS_\rho$, up 
to translation we can always assume that $y=0 \notin \cup_i B^2_i$, where
 $\rho\ge 2$ is a fixed parameter that will be chosen large enough at some point in the proof.
Throughout the proof we impose the relation $r:=\rho/2$.
We also denote by $C_d\gg 1$ a universal constant (that may change from line to line but can be chosen depending on $d$ only).

\medskip

\step1 Structure of the proof.

We start by writing the PDE solved by $x\mapsto G(x)$ on $\R^d$, namely
$$
(\frac1T-\triangle) G(x,y) \,=\,\delta (x)+\frac1T G \mathds 1_{\calB} - \sum_i \partial_n G \delta_{\partial B_i}.
$$
To avoid dealing with the singularity of the RHS, we convolve $G$ with the simple moving average $m$ on the unit all $B$, and set $g(x)=\fint_B G(x+z)dz=m*G$.
The averaged function $g$ then satisfies
\begin{equation}\label{e-green-1.1}
(\frac1T-\triangle) g(x) \,=\, m(x)+\frac1T m*(G \mathds 1_{\calB}) - \sum_i m*( \partial_n G\delta_{\partial B_i}).
\end{equation}
We finally consider the standard massive Green's function $h$ solution on $\R^d$ of 
\begin{equation}\label{e-green-1.2}
(\frac1T-\triangle) h(x) \,=\,\delta (x),
\end{equation}
and that satisfies for $d>2$ the estimates 
\begin{equation}\label{e-green-1.5}
\sum_{j=0}^3 |x|^j |\nabla^j h(x)| \lesssim |x|^{2-d} e^{-\frac{|x|}{C_d \sqrt T}}.
\end{equation}
The combination of \eqref{e-green-1.1} and \eqref{e-green-1.2} then yields
\begin{eqnarray}\nonumber
g(x) &=& m*h(x)+\frac1T \int G \mathds 1_{\calB} (m*h)(x-\cdot)\\
&&- \sum_i \int_{\partial B_i} (\partial_n G) (m*h)(x-\cdot), \label{e-green-1.3}
\\
\nabla g(x) &=& \nabla m*h(x)+\frac1T \int G \mathds 1_{\calB} \nabla (m*h)(x-\cdot)\nonumber
\\
&&-\sum_i \int_{\partial B_i} (\partial_n G) \nabla (m*h)(x-\cdot).\label{e-green-1.4}
\end{eqnarray}
Since $\int_{\partial B_i}\partial_n G=0$, $m*h$ and $\nabla (m*h)$ in the RHS sums of \eqref{e-green-1.3}
and \eqref{e-green-1.4} can be replaced by  $\nabla (m*h)$ and  $\nabla^2 (m*h)$, respectively, when it turns to 
estimates. However,  $\nabla^2 (m*h)$ is borderline \emph{non-integrable}, which makes the analysis more subtle and requires
to unravel further cancellations.
The aim of the upcoming steps is to prove that the RHS terms of \eqref{e-green-1.3} \& \eqref{e-green-1.4} that involve $G$ can be absorbed into the LHS by means of
Neumann series. 
More precisely we shall focus on the quantities
$$g_i:=|G_i| \quad \text{and} \quad  \gamma_i:=\Big(\int_{B_i^r} |\nabla G|^2\Big)^\frac12,$$
and reformulate \eqref{e-green-1.3} and \eqref{e-green-1.4} to obtain iterable estimates
at the level of $g_i$ and $\gamma_i$.

\medskip

\step2 Estimate of $\int_{\partial B_i} (\partial_n G) (m*h)(x-\cdot)$ and $\int_{\partial B_i} (\partial_n G) \nabla (m*h)(x-\cdot)$.

\nopagebreak

Let $x \in B^r_j$ for $j$ fixed.
We first consider $i\ne j$.
Since $m*h$ is smooth we have by \eqref{e-green-1.5} 
 (using that $|x_i-x_j|\gtrsim \inf_{z\in \partial B_i} |x-z|$ since $\rho=2r$, and $|x-x_j|\le r$):
 For all $z\in \partial B_i$,
\begin{eqnarray*}
&&{|m*h(x-z)-m*h(x-x_i))| }
\,\lesssim \, |x_j-x_i|^{1-d}e^{-\frac{|x_j-x_i|}{C_d \sqrt T}},
\\
&&{|\nabla m*h(x-z)-\nabla m*h(x-x_i)-(\nabla m*h(x_j-z)-\nabla m*h(x_j-x_i))| }
 \\
 && \hspace{5cm}\,\lesssim \, \rho |x_j-x_i|^{-d-1}e^{-\frac{|x_j-x_i|}{C_d \sqrt T}}.
\end{eqnarray*}
Using first $\int_{\partial B_i} \partial_n G_T=0$, 
this yields
\begin{eqnarray*}
|\int_{\partial B_i} (\partial_n G) (m*h)(x-\cdot)| &\lesssim & |x_j-x_i|^{1-d} e^{-\frac{|x_j-x_i|}{C_d \sqrt T}}\int_{\partial B_i} |\partial_n G|, 
\\
|\int_{\partial B_i} (\partial_n G) \Big( \nabla (m*h)(x-\cdot)-\nabla (m*h)(x_j-\cdot)\Big)|  &\lesssim &\rho |x_j-x_i|^{-d-1}e^{-\frac{|x-x_i|}{C_d \sqrt T}} \int_{\partial B_i} |\partial_n G|.
\end{eqnarray*}
We then appeal to Lemma~\ref{lem:reg1} and obtain in the regime $\sqrt T \ge \rho=2r$
\begin{eqnarray}
\lefteqn{|\int_{\partial B_i} (\partial_n G) (m*h)(x-\cdot)| }\nonumber
\\
&\lesssim & \rho^{-\alpha} |x_j-x_i|^{1-d} e^{-\frac{|x_j-x_i|}{C_d \sqrt T}}\Big(\int_{B_i^r} |\nabla G|^2+ \rho^d \frac{G_i^2}{T^2}\Big)^\frac12, \label{e-green-2.1}
\\
\lefteqn{|\int_{\partial B_i} (\partial_n G) \Big( \nabla (m*h)(x-\cdot)-\nabla (m*h)(x_j-\cdot)\Big)|  }\nonumber
\\
&\lesssim &\rho^{1-\alpha} |x_j-x_i|^{-d-1}e^{-\frac{|x_j-x_i|}{C_d \sqrt T}} \Big(\int_{B_i^r} |\nabla G|^2+ \rho^d \frac{G_i^2}{T^2}\Big)^\frac12.
\label{e-green-2.2}
\end{eqnarray}
For $i=j$, we do not replace $x$ by $x_j$ and rather obtain by a similar string of arguments as above
\begin{eqnarray}
\lefteqn{|\int_{\partial B_i} (\partial_n G) (m*h)(x-\cdot)| }\nonumber
\\
&\lesssim & \rho^{-\alpha} (1+|x-x_i|)^{1-d} e^{-\frac{|x-x_i|}{C_d \sqrt T}}\Big(\int_{B_i^r} |\nabla G|^2+ \rho^d \frac{G_i^2}{T^2}\Big)^\frac12, \label{e-green-2.1bis}
\\
\lefteqn{|\int_{\partial B_i} (\partial_n G) \nabla (m*h)(x-\cdot)|  }\nonumber
\\
&\lesssim &\rho^{-\alpha} (1+|x-x_i|)^{-d}e^{-\frac{|x-x_i|}{C_d \sqrt T}} \Big(\int_{B_i^r} |\nabla G|^2+ \rho^d \frac{G_i^2}{T^2}\Big)^\frac12.
\label{e-green-2.2bis}
\end{eqnarray}
(The scalings are different in \eqref{e-green-2.2} and \eqref{e-green-2.2bis}.)
We also note that for $x=x_j$, we have incidentally the neat identity
\begin{equation}\label{e-green-2.3}
\int_{\partial B_j} (\partial_n G) (m*h)(x_j-\cdot)=0
\end{equation}
since $h$ (and therefore $m*h$) is radially symmetric, $B_j$ is centered at $x_j$,  and $\int_{\partial B_j} \partial_n G=0$.
\medskip

\step3 Iterable estimates.

Since $G$ is constant on balls $B_j$, $g(x_j)=\fint_{B_j} G=G_j$.
The combination of \eqref{e-green-1.3}, \eqref{e-green-1.5}, \eqref{e-green-2.1}, and \eqref{e-green-2.3} yields for all $x_j$
the nonlinear estimate
\begin{multline}\label{e-green-4.1bis}
|G_j|=|g(x_j)| \,\lesssim \, |x_j|^{2-d} e^{-\frac{|x_j|}{C_d \sqrt T}}+ \sum_{i}  |G_i| \frac1T(1+|x_j-x_i|)^{2-d} e^{-\frac{|x_j-x_i|}{C_d \sqrt T}}
\\
 + \rho^{-\alpha} \sum_{i\ne j} |x_j-x_i|^{1-d} e^{-\frac{|x_j-x_i|}{C_d \sqrt T}}\Big(\int_{B_i^r} |\nabla G|^2+\rho^d \frac{G_i^2}{T^2}\Big)^\frac12.
\end{multline}
We then turn to the gradient. Let $x \in B^r_j$.
Since $G$ is constant on $B_j$, 
$$
0=\nabla g(x_j)=
\nabla m*h(x_j)+\frac1T \int G \mathds 1_{\calB} \nabla (m*h)(x_j-\cdot)-\sum_i \int_{\partial B_i} (\partial_n G) \nabla (m*h)(x_j-\cdot),
$$
and therefore by the triangle inequality
\begin{multline*}
|\nabla g(x)| = |\nabla g(x)-\nabla g(x_j)|
\\\lesssim \, |\nabla m*h(x)| + |\nabla m*h(x_j)| +\frac1T \int G \mathds 1_{\calB}( |\nabla (m*h)(x-\cdot)|+ |\nabla (m*h)(x_j-\cdot)|)
\\
+\sum_i \Big|\int_{\partial B_i} (\partial_n G) \Big(\nabla (m*h)(x-\cdot)-\nabla (m*h)(x_j-\cdot)\Big)\Big|.
\end{multline*}
Combined with \eqref{e-green-1.4}, \eqref{e-green-1.5}, \eqref{e-green-2.2}, and \eqref{e-green-2.2bis}, this entails 
\begin{multline}\label{e-green-4.2}
|\nabla g(x)| \,\lesssim \, |x_j|^{1-d} e^{-\frac{|x_j|}{C_d \sqrt T}}+ |G_j| \frac1T(1+|x-x_j|)^{1-d}e^{-\frac{|x-x_i|}{C_d \sqrt T}}
+ \sum_{i\ne j} |G_i| \frac1T |x_j-x_i|^{1-d}e^{-\frac{|x_j-x_i|}{C_d \sqrt T}}
\\
 +\rho^{-\alpha}  (1+|x-x_j|)^{-d} e^{-\frac{|x-x_i|}{C_d \sqrt T}}\Big(\int_{B_j^r} |\nabla G|^2+\rho^d \frac{G_i^2}{T^2}\Big)^\frac12
 \\
+ \rho^{1-\alpha}  \sum_{i\ne j} |x_j-x_i|^{-d-1} e^{-\frac{|x_j-x_i|}{C_d \sqrt T}}\Big(\int_{B_i^r} |\nabla G|^2+\rho^d \frac{G_i^2}{T^2}\Big)^\frac12.
\end{multline}
We now reformulate \eqref{e-green-4.1bis} and \eqref{e-green-4.2} into iterable estimates at the level 
of the quantities $G_j$ and $\gamma_j$.
Estimate  \eqref{e-green-4.1bis} directly takes the form 
\begin{multline}\label{e-green-4.3}
g_j \,\lesssim \, |x_j|^{2-d} e^{-\frac{|x_j|}{C_d \sqrt T}}+ \sum_i  \frac1T(1+|x_j-x_i|)^{2-d}e^{-\frac{|x_j-x_i|}{C_d \sqrt T}}g_i 
\\
 + \rho^{-\alpha} \sum_{i\ne j} |x_j-x_i|^{1-d} e^{-\frac{|x_j-x_i|}{C_d \sqrt T}} (\gamma_i+\frac{\rho^\frac d2}{ T} g_i).
\end{multline}
We then turn to \eqref{e-green-4.2} and let $x \in B^r_j \setminus B^2_j$.
By the mean value property for the massive Laplacian (see e.~g.~\cite[(39) page 289]{Courant-Hilbert-II}),
\begin{multline}\label{e.a+1}
|\nabla G(x)|\, \lesssim \, |\nabla g(x)| \,\lesssim \, |x_j|^{1-d} e^{-\frac{|x_j|}{C_d \sqrt T}}+ 
 \frac1T(1+|x-x_j|)^{1-d}e^{-\frac{|x-x_j|}{C_d \sqrt T}}g_j  
 \\
 +\sum_{i\ne j}  \frac1T|x_j-x_i|^{1-d}e^{-\frac{|x_j-x_i|}{C_d \sqrt T}}g_i 
\\
+\rho^{-\alpha}  (1+|x-x_j|)^{-d} e^{-\frac{|x-x_i|}{C_d \sqrt T}} (\gamma_j+\frac{\rho^\frac d2}{ T} g_j)
 \\
 + \rho^{1-\alpha} \sum_{i\ne j} |x_j-x_i|^{-d-1} e^{-\frac{|x_j-x_i|}{C_d \sqrt T}} (\gamma_i+\frac{\rho^\frac d2}{ T} g_i).
\end{multline}
We now reconstruct $\gamma_j$ using $|\nabla G(x)|$ for $x\in B_j^r \setminus B_j^2$ and proceed in two steps.
First we integrate the square of \eqref{e.a+1} over $B^r_j\setminus B^2_j$ to the effect of
\begin{multline*}
\int_{B^r_j\setminus B^2_j} |\nabla G(x)|^2dx \, \,\lesssim \, \rho^d \bigg(|x_j|^{1-d} e^{-\frac{|x_j|}{C_d \sqrt T}}+ \sum_{i \ne j}  \frac1T|x_j-x_i|^{1-d}e^{-\frac{|x_j-x_i|}{C_d \sqrt T}}g_i 
\\
 + \rho^{1-\alpha}  \sum_{i\ne j} |x_j-x_i|^{-d-1} e^{-\frac{|x_j-x_i|}{C_d \sqrt T}} (\gamma_i+\frac{\rho^\frac d2}{ T} g_i)\bigg)^2+ \rho^{-2\alpha}(\gamma_j+\frac{\rho^{\frac d2}}{T} g_j)^2.
\end{multline*}
Second, we appeal to the energy estimate (which we prove at the end of this step) 
\begin{equation}\label{e-green-4.4}
\int_{B^2_j\setminus B_j} |\nabla G(x)|^2dx\,\lesssim\, \int_{B^3_j\setminus B^2_j} |\nabla G(x)|^2dx + \frac1T g_j^2.
\end{equation}
The combination of these two estimates with the hardcore condition $|x_i-x_j|\ge \rho \delta(i-j)$ 
in form of $\rho^{1-\alpha}   |x_j-x_i|^{-d-1}\le \rho^{-\frac \alpha2 } |x_j-x_i|^{-d- \frac \alpha2}$
finally yields the iterable estimate 
\begin{multline}\label{e-green-4.5}
\gamma_j \,\lesssim \, \rho^\frac d2 |x_j|^{1-d} e^{-\frac{|x_j|}{C_d \sqrt T}}+  \sum_{i\ne j}  \frac{\sqrt \rho}T (\frac{|x_j-x_i|}{\sqrt \rho})^{1-d}e^{-\frac{|x_j-x_i|}{C_d \sqrt T}}g_i 
\\
 + \rho^{-\frac \alpha 4}   \sum_{i \ne j} (\frac{|x_j-x_i|}{\sqrt \rho})^{-d-\frac \alpha2} e^{-\frac{|x_j-x_i|}{C_d \sqrt T}} (\gamma_i+\frac{\rho^\frac d2}{ T}g_i)+ \rho^{-\alpha}(\gamma_i+\frac{\rho^\frac d2}{ T} g_i)+ \frac1{\sqrt T} g_j.
\end{multline}
We conclude this step with the argument in favor of \eqref{e-green-4.4}.
We recall the following two properties of the $(\frac1T-\triangle)$-harmonic function $G$ on $B^2_j \setminus B_j$: For all $2\le t \le 3$,
since $\int_{\partial B_j} \partial_n G=0$ and $G\equiv G_j$ on $\partial B_j$,
\begin{equation}\label{e.a+2}
\int_{B_j^t \setminus B_j} \frac1T G -\int_{\partial B_j^t} \partial_n G+\int_{\partial B_j} \partial_n G=0 \quad \implies \quad \int_{\partial B_j^t} \partial_n G=\int_{B_j^t \setminus B_j} \frac1T G,
\end{equation}
and 
\begin{multline}\label{e.a+3}
\int_{B^t_j \setminus B_j} \frac1T G^2+\int_{B_j^t \setminus B_j} |\nabla G|^2 -\int_{\partial B_j^t} G\partial_n G+\int_{\partial B_j}G \partial_n G=0\quad 
\\
\implies \quad \int_{B_j^t \setminus B_j} |\nabla G|^2  \le \int_{\partial B_j^t} G\partial_n G .
\end{multline}
The combination of these two properties entails by Poincar\'e's inequality on $\partial B_j^t$ for $2\le t\le 3$
and on $B^3_j\setminus B_j$
\begin{eqnarray*}
\int_{B^2_j\setminus B_j}|\nabla G|^2  &\le& \int_2^3  \int_{B^t_j \setminus B_j} |\nabla G|^2 dt 
\\
& \stackrel{\eqref{e.a+2}\& \eqref{e.a+3}}\le& \int_2^3 \int_{\partial B_j^t} (G-\fint_{\partial B_j^t} G)\partial_n Gdt 
+\int_2^3 (\fint_{\partial B_j^t} G) \int_{B_j^t \setminus B_j} \frac1T G dt
\\
&\lesssim & \int_2^3 \int_{\partial B_j^t} |\nabla G|^2 dt+  \int_2^3 \int_{\partial B_j^t} |G| \int_{B_3 \setminus B_j} \frac1T |G| dt
\\
&\lesssim& \int_{B^3_j\setminus B^2_j}  |\nabla G|^2 +\frac1T \int_{B^3_j} G^2
\\
&\lesssim &   \int_{B^3_j\setminus B^2_j}  |\nabla G|^2 +\frac1T g_j^2 + \frac1T \int_{B^3_j \setminus B_j} |\nabla G|^2.
\end{eqnarray*}
The desired estimate \eqref{e-green-4.4} then follows from absorbing part of the last RHS term into the LHS for $T\gg 1$.

\medskip

\step4 Estimates of $g_i$ and $\gamma_i$.

\nopagebreak
In view of the hardcore condition on the point process, we may wlog parametrize the point set by  $\Z^d$ (that is, a point $x \in \calP$
is labelled $x_i$ for $i \in \arg\inf\{ |x- c_d \rho j| \, |\,  j \in  \Z^d \}$, which  defines an injective --- although not surjective --- enumeration for $c_d$ small enough).
Set $\Gamma_j= \frac{\rho^\frac d2}{\sqrt T}g_j+\gamma_j$ if $x_j \in \calP$ and $\Gamma_j=0$ otherwise, and set $B_j = \rho^{\frac d2}(1+\rho |j|)^{1-d}e^{-\frac{\rho |j|}{C_d\sqrt T}}$ for all $j \in \Z^d$.
The combination of  \eqref{e-green-4.3} and \eqref{e-green-4.5} then yields
\begin{equation}\label{e-green-5.1}
\Gamma_j \,\lesssim\, B_j+M_{jj} \Gamma_j+\sum_{i \ne j} M_{ij} \Gamma_i,
\end{equation}
where 
\begin{eqnarray*}
M_{jj}&=&(\frac1T+\rho^{-\alpha}+\frac1{\sqrt \rho^{d}}), \\
M_{ij}&=& \bigg( \frac{1}{\sqrt{\rho}^{d}} \frac \rho T(\sqrt \rho |j-i|)^{2-d} + \rho^{-\alpha}\frac{\sqrt \rho}{\sqrt T}(\sqrt \rho {|j-i|})^{1-d}
\\
&&\hspace{2cm} + \frac{1}{\sqrt \rho^d}\frac{\sqrt \rho}{\sqrt T} (\sqrt \rho {|j-i|})^{1-d}+\rho^{-\frac \alpha4}  (\sqrt \rho |j-i|)^{-d-\frac \alpha2}\bigg)e^{-\frac{\sqrt \rho |j-i|}{C_d \sqrt T/\sqrt \rho}}.
\end{eqnarray*}
It is now easy to check that for all $C\gg 1$,  provided   $\sqrt T\ge  \rho \gg_C 1$, we have both
$$M_{ij} \le \frac1C (1+\sqrt \rho |i-j|)^{-d} e^{-\frac{\sqrt \rho |j-i|}{C_d \sqrt T/\sqrt \rho}} \quad \text{ and } \quad \sup_j \sum_i M_{ij} \le \frac1C$$
(note that the first estimate does not imply the second one).
For $(D)_j$ given by  $D_j= \rho^\frac d2 (1+\sqrt \rho |j|)^{1-d}e^{-\frac{\sqrt \rho |j|}{C_d\sqrt T/\sqrt \rho}} \gtrsim B_j$ for all $j$, this implies (for some suitable choice of $C\gg 1$)
$$
(M D)_j \le \frac12 D_j.
$$
Therefore $(\sum_{k=0}^\infty M^k D)_j \le 2 D_j$, and thus $(\sum_{k=0}^\infty M^k B)_j \le 2 D_j$.
Hence, \eqref{e-green-5.1} entails for all $N\in \N$
\begin{equation}\label{e-green-5.3}
\Gamma \,\lesssim\, (\sum_{k=0}^\infty M^k D) +M^N \Gamma \,\lesssim\, D+M^N \Gamma.
\end{equation}
By a variant of the energy estimates of Lemma~\ref{lem:massive} (that we display for completeness
in Step~2 of the proof of Lemma~\ref{lem:massive} below), one has for all $y \in \R^d \setminus \cup_{i}B^2_i$ and all $|x-y|\ge 2$
the (largely suboptimal) estimate
\begin{equation}\label{e.Green-bebe}
\frac1{\sqrt T} |G(x,y)|+|\nabla G(x,y)| \,\lesssim\, \sqrt{T} e^{-\frac{|x-y|}{C_d\sqrt T}},
\end{equation}
which implies $\Gamma_j \lesssim \rho^\frac d2\sqrt{T}^{d}(\rho |j|+1)^{1-d} \exp(-\frac{\rho |j|} {C_d \sqrt T})$.
This coarse estimate is however enough to pass to the limit $N\uparrow \infty$ in \eqref{e-green-5.3}, which yields the crucial estimate
\begin{equation}\label{e-green-5.2}
\frac{\rho^\frac d2}{\sqrt T}g_j+\gamma_j= \Gamma_j  \, \lesssim \, D_j  \lesssim \, \rho^\frac d2 (1+\frac{|x_j|}{\sqrt \rho})^{1-d}e^{-\frac{|x_j|}{C_d\sqrt T}} .
\end{equation}

\medskip

\step5 Conclusion: Proof of \eqref{e.Green1} and \eqref{e.Green2}.

Since for all $j$,
$$
\int_{\partial B_j} |\partial_n G|\,\lesssim\, \gamma_j,
$$
the estimate
$$
|g(x)|\,\lesssim\,   (1+\frac{|x|}{\sqrt \rho})^{2-d}e^{-\frac{|x|}{C_d\sqrt T}} 
$$
follows for all $|x|\ge   2$ by \eqref{e-green-5.2}, \eqref{e-green-1.3}, and the property $\int_{\partial B_i} \partial_n G=0$.
To get rid of the moving average, and deduce \eqref{e.Green1}, it suffices to appeal to the mean-value property (for the massive Laplacian) away from the inclusions, and to local elliptic regularity on the remaining regions.
\medskip

We conclude with the argument in favor of \eqref{e.Green2}.
If $x\in B_j^{r-1}$, this follows from \eqref{e-green-5.2} by the mean-value property (for the massive Laplacian) away from the inclusion $B_j$, and by elliptic regularity on the remaining regions.
It remains to argue in the case when $x \notin \cup_{j}B^{r-1}_j$. 
We further distinguish the case $|x|\le \frac12 \dist(0,\calP_\rho)=:s$
and $|x|\ge s$.
If $2\le |x|\le s$, we obtain by Caccioppoli's estimate on $B_{|x|/2}(x)$ (on which $G$ is $(\frac1T-\triangle)$-harmonic):
$$
\fint_{B_{|x|/4}(x)} |\nabla G|^2 \,\lesssim\, (|x|^{-2}+\frac1T) \fint_{B_{|x|/2}(x)}G^2  \,\stackrel{\eqref{e.Green1}}\lesssim\, |x|^{2(1-d)} e^{-\frac{|x|}{C_d \sqrt T}},
$$
and the claim follows from the mean-value property. 
Let $1\le k \le d$.
Combined with \eqref{e-green-5.2} and the explicit formula for the Green's function of the massive Laplacian, the above implies that there exist $T' \ge T$ and $C<\infty$ such that for all $x\in \cup_{j}B^{r-1}_j \cup B^s(0)$,
$|\partial_k G(x)| \le C \sum_j |\partial_j h_{T'}(x)|$.
Since 
$|\partial_j h_{T'}|$ is $\frac1{T'}-\triangle$-superharmonic on $\R^d \setminus (\cup_i B^{r-1}_i \cup B^s(0))$ and $\partial_k G$ is $\frac1T-\triangle$-harmonic on $\R^d \setminus (\cup_i B^{r-1}_i \cup B^s(0))$, we obtain
$$
\frac1T(\partial_k G-C\sum_j|\partial_j h_{T'}|)-\triangle (\partial_i G-C\sum_j|\partial_i h_{T'}|) \le C( \frac{1}{T'}-\frac1{T})\sum_j|\partial_j h_{T'}|\le 0,
$$
so that $\partial_k G\le C\sum_j |\partial_j h_{T'}|$ by the maximum principle.
Likewise $-\partial_k G\le C\sum_j |\partial_j h_{T'}|$, and therefore $|\partial_k G|\le C\sum_j |\partial_j h_{T'}|$. The desired estimate \eqref{e.Green2} follows.

\section{Proofs of the other auxiliary results}\label{sec:auxiliary-results}

\subsection{Lemma~\ref{lem:massive}: Deterministic existence result}

For all $C>0$, we consider the exponential cut-off $\tilde \eta_{T,C}(x):=\exp(-\frac{|x|}{C\sqrt{T}})$; it satisfies $|\nabla \tilde \eta_{T,C}|\lesssim \frac1{C\sqrt{T}}\tilde \eta_{T,C}$. 
We shall modify this cut-off so that it be constant on the balls $B_i$.
For all $i$, define open balls $\tilde B_i$ and $\bar B_i$ centered at $x_i$s in such a way that $B_i \subset\subset \tilde B_i \subset \subset \bar B_i$ and
that the balls $\bar B_i$ are all disjoint, and set $\bar \calB=\cup_i \bar B_i$ and $\R^d_{\bar \calB}=\R^d\setminus \bar \calB$.
Define $\eta_{T,C}|_{\R^d_{\bar \calB}}$ as follows:
\begin{itemize}
\item $\eta_{T,C}|_{\R^d_{\bar \calB}} \equiv \tilde \eta_{T,C}|_{\R^d_{\bar \calB}}$, 
\item For all $i$: $\eta_{T,C}|_{\tilde B_i}\equiv \inf_{\bar B_i} \tilde \eta_{T,C}$,
\item For all $i$: Extend $\eta_{T,C}$ linearly radially with respect to $x_i$ between $\partial \bar B_i$ and $\partial \tilde B_i$. 
\end{itemize}
So defined,  for all $v\in \Hc_\uloc$, $\eta_{T,C}^2 v \in \Hc_\uloc$, whereas we still have the crucial estimate 
$|\nabla \eta_{T,C}|\lesssim \frac1{C\sqrt{T}} \eta_{T,C}$ on $\R^d$.
 
\medskip
 
\step1 Caccioppoli argument and proof of  \eqref{e.massive-apriori}.
 
We are in the position to proceed to the Caccioppoli argument. We test the equation \eqref{e.corr-Tdet} with $\eta_{T,C}^2 v_T$ and integrate on $\R^d$.
This yields after integrations by parts
$$
\frac1T\int_{\R^d\setminus\calB} \eta_{T,C}^2 v_T^2 +\int_{\R^d\setminus\calB} |\nabla (\eta_{T,C}v_T)|^2 = \sum_i \eta_{T,C,i}^2 g_2 v_{T,i} +\int_{\R^d\setminus\calB} |\nabla \eta_{T,C}|^2 v_T^2 +\int_{\R^d\setminus\calB}  g_1 \eta_{T,C}^2 v_T,
$$
where $\eta_{T,C,i}$ denotes the value of the cut-off on $B_i$.
Let $C'>0$ denote a finite constant that may change from line to line, and will be chosen to be large enough at the end of the proof.
For the first RHS we use the definition of $\eta_{T,C}$ and a trace estimate on each $\partial B_i$ in form of 
\begin{eqnarray*}
\Big|\sum_i \eta_{T,C,i}^2  g_2 v_{T,i}\Big|&\lesssim &  \sum_i \Big(\int_{\bar B_i \setminus B_i}\eta_{T,C}^2  g_2^2 \Big)\Big(\int_{\bar B_i \setminus B_i} (\eta_{T,C} v_T)^2+|\nabla (\eta_{T,C} v_T)|^2\Big)^\frac12
\\
&\lesssim  & C'T g_2^2 \int_{\bar \calB \setminus \calB} \eta_{T,C}^2+ \frac1{C'T}  \int_{\R^d\setminus\calB}\eta_{T,C}^2 v_T^2+|\nabla (\eta_{T,C} v_T)|^2.
\end{eqnarray*}
Likewise, for the last RHS term
$$
\Big|\int_{\R^d\setminus\calB}  g_1 \eta_{T,C}^2 v_T\Big| \,\lesssim\, \frac{1}{C'T} \int_{\R^d\setminus\calB} \eta_{T,C}^2 v_T^2+C'T g_1^2\int_{\R^d\setminus\calB} \eta_{T,C}^2,
$$
whereas for the second RHS term we use the property of the cut-off
$$
\int_{\R^d\setminus\calB} |\nabla \eta_{T,C}|^2 v_T^2  \,\lesssim\, \frac1{C^2T} \int_{\R^d\setminus\calB}  \eta_{T,C}^2 v_T^2.
$$
Choosing $C$ and $C'$ large enough to absorb part of the RHS into the LHS, these last four estimates combine to 
\begin{equation}\label{e.energy-massive-refined}
\frac1T\int_{\R^d\setminus\calB} \eta_{T,C}^2 v_T^2 +\int_{\R^d\setminus\calB} |\nabla (\eta_{T,C}v_T)|^2 \,\lesssim\,C'T g_2^2 \int_{\bar \calB \setminus \calB} \eta_{T,C}^2
+C'T g_1^2\int_{\R^d\setminus\calB} \eta_{T,C}^2,
\end{equation}
which implies \eqref{e.massive-apriori}  since the origin plays no role in this estimate.
This yields existence and uniqueness by standard arguments (solutions can be constructed by approximation on balls of radius $r>0$ with 
homogeneous Dirichlet boundary conditions, and the uniform estimate allows to take the limit $r\uparrow \infty$).

\medskip

We conclude this proof by an auxiliary result that is needed as a mild starting point to establish the sharp decay of the Green's function $G_T$.

\medskip

\step2 Suboptimal a priori estimates for the Green's function.

Let $y \in \R^d \setminus \cup_{B_i^2}$, and consider $x\mapsto G(x,y)$ the solution of
\begin{equation*} 
\frac1TG -\triangle G\,=\, \mathds 1_{B(y)} \quad \text{ in } \R^d \setminus \calB, \quad  \forall i:\quad G|_{B_i} \equiv G_{i} \in \R,  \quad \int_{\partial B_i} \partial_n G = 0.
\end{equation*}
The function $G$ is a moving average of the Green's function (for which $\mathds 1_{B(y)}$ is replaced by a Dirac mass at $y$).
Wlog we assume $x=0$.
The argument leading to \eqref{e.energy-massive-refined} yields in this case
$$
\frac1T\int_{\R^d\setminus\calB} \eta_{T,C}^2 G^2 +\int_{\R^d\setminus\calB} |\nabla (\eta_{T,C}G)|^2 \,\lesssim\,C'T  \int_{B(y)} \eta_{T,C}^2.
$$
This entails 
$$
\int_{B(0)}\frac1T G ^2 + |\nabla G|^2 \,\lesssim\, T e^{-\frac{|y|}{C_d \sqrt T}}.
$$
To get rid of the moving average, it suffices to appeal to the mean-value property (for the massive Laplacian) away from the inclusions, and to energy estimates on the remaining regions.
The argument holds for $0$ replaced by any $x$ with $|x-y|\ge 2$. This proves \eqref{e.Green-bebe}.

\subsection{Lemma~\ref{lem:massive-random}: Energy estimates for the massive corrector}

By Lem\-ma~\ref{lem:massive}, we have existence and uniqueness of $u_T$ almost surely. 
Uniqueness implies that $u_T$ is stationary, as claimed.
Both \eqref{e.massive-apriori-random} and \eqref{e.expec-V_iT} follow from the weak form \eqref{e.eq-weak-prob} of the equation in probability.
Indeed, taking $v\equiv 1$ yields
\begin{equation*}
\frac1T \expec{u_T \mathds1_{\R^d\setminus\calB}} \,=\, \expec{\bar g_\theta \mathds1_{\R^d\setminus\calB}-\bar g \mathds1_{\calB}}=(1-\theta) \bar g_\theta-\theta \bar g = 0
\end{equation*}
by the choice of $g_\theta$, that is, \eqref{e.eq-expec-uT}.
By taking $v=u_T$ in \eqref{e.eq-weak-prob}, we obtain \eqref{e.expec-V_iT} using \eqref{e.eq-expec-uT}.
By a trace estimate with constant $C''>0$ and Young's inequality
\begin{eqnarray*}
 |\bar g\expec{u_T\mathds 1_{\calB}}|&\le& |\bar g| \expec{\mathds 1_{\calB}}^\frac12 \expec{u_T^2\mathds 1_{\calB}}^\frac12
 \\
 &\le& \frac{C'}2 T\theta \bar g^2+\frac{C''}{2TC'} \expec{(u_T^2+|\nabla u_T|^2) \mathds 1_{\R^d \setminus \calB}},
\end{eqnarray*}
where $C'$ is arbitrary and chosen so small that the second RHS term can be absorbed in the LHS of \eqref{e.expec-V_iT} to the effect that
$$
\expec{(\frac1T u_T^2+|\nabla u_T|^2)\mathds1_{\R^d\setminus\calB}} \lesssim  T \theta \bar g^2 .
$$
Estimate~\eqref{e.massive-apriori-random} then follows from the observation that $\nabla u_T=0$ on $\calB$ and from the trace estimate.
It remains to establish \eqref{e.eq-weak-prob}.
Consider the cut-off $\eta_{T,C}$ of the proof of Lemma~\ref{lem:massive}, and test \eqref{e.corr-T} with function $\eta_{T,C} v$.
This yields after integration by parts 
\begin{multline}\label{e.a+4}
\frac1T \int_{\R^d}  \eta_{T,C}u_T v \mathds1_{\R^d\setminus \calB}+\int_{\R^d} \eta_{T,C} \nabla u_T \cdot \nabla v \,=\, \int_{\R^d} \eta_{T,C}v(  \bar g_\theta \mathds 1_{\R^d\setminus\calB}-\bar g \mathds 1_{\calB} )
\\
-\int_{\R^d} v \nabla \eta_{T,C} \cdot \nabla u_T.
\end{multline}
We first note that for all $x\in \R^d$, $|\eta_{T,C}(x)-\tilde \eta_{T,C}(x)|\lesssim |\nabla \tilde \eta_{T,C}(x)|\lesssim \frac1{C\sqrt{T}} \tilde \eta_{T,C}(x)$. 
Hence, by taking expectations, we have 
\begin{multline*}
\Big|\expec{\frac1T \int_{\R^d}  (\eta_{T,C}-\tilde \eta_{T,C})u_T v \mathds1_{\R^d\setminus \calB}+\int_{\R^d} (\eta_{T,C}-\tilde \eta_{T,C}) \nabla u_T \cdot \nabla v}
\\
-\expec{ \int_{\R^d} (\eta_{T,C}-\tilde \eta_{T,C})v(\bar g_\theta \mathds 1_{\R^d\setminus\calB}- \bar g \mathds 1_{\calB}) }\Big|
\\
\lesssim \, 
\frac1{C\sqrt{T}} \expec{ \int_{\R^d}  \tilde \eta_{T,C}( \frac1Tu_T^2+  \frac1T|v|^2+|\nabla u_T|^2+ |\nabla v|^2) + \int_{\R^d}  \tilde \eta_{T,C} |v| |\bar g| }
\\
+\frac1{C\sqrt{T}} \expec{\int_{\R^d} \tilde \eta_{T,C} |v| |\nabla u_T| }.
\end{multline*}
On the one hand, since $v$ and $u_T$ are stationary and $\tilde \eta_{T,C}$ is deterministic, this yields
\begin{multline*}
\Big|\expec{\frac1T \int_{\R^d}  (\eta_{T,C}-\tilde \eta_{T,C})u_T v \mathds1_{\R^d\setminus \calB}+\int_{\R^d} (\eta_{T,C}-\tilde \eta_{T,C}) \nabla u_T \cdot \nabla v}
\\
-\expec{ \int_{\R^d} (\eta_{T,C}-\tilde \eta_{T,C})v(\bar g_\theta \mathds 1_{\R^d\setminus\calB}- \bar g \mathds 1_{\calB} )}\Big|
\\
\lesssim \, 
\frac1{C\sqrt{T}} \expec{\frac1Tu_T^2+  v^2+\bar g^2+|\nabla u_T|^2+ |\nabla v|^2 } \int_{\R^d}  \tilde \eta_{T,C}.
\end{multline*}
On the other hand, again by stationarity and the fact that $\tilde \eta_{T,C}$ is deterministic (and in contrast to $\eta_{T,C}$), 
\begin{multline*}
\expec{\frac1T \int_{\R^d} \tilde \eta_{T,C} u_T v \mathds1_{\R^d\setminus \calB}+\int_{\R^d}\tilde \eta_{T,C} \nabla u_T \cdot \nabla v
- \int_{\R^d}\tilde \eta_{T,C}v(  \bar g_\theta \mathds 1_{\R^d\setminus\calB}-\bar g \mathds 1_{\calB} )}
\\
= \expec{\frac1Tu_T v \mathds1_{\R^d\setminus \calB}+ \nabla u_T \cdot \nabla v
-v( \bar g_\theta \mathds 1_{\R^d\setminus\calB}-\bar g \mathds 1_{\calB})} \int_{\R^d} \tilde \eta_{T,C} ,
\end{multline*}
and by \eqref{e.a+4} and since $|\nabla \eta_{T,C}|\lesssim \frac{1}{C\sqrt T} \tilde \eta_{T,C}$ is a deterministic bound,
\begin{multline*}
\Big|\expec{\frac1T \int_{\R^d}  \eta_{T,C} u_T v \mathds1_{\R^d\setminus \calB}+\int_{\R^d} \eta_{T,C} \nabla u_T \cdot \nabla v
- \int_{\R^d} \eta_{T,C}v( \bar g \mathds 1_{\calB} - \bar g_\theta \mathds 1_{\R^d\setminus\calB})}\Big|
\\
\lesssim \, \frac1{C\sqrt{T}} \expec{ v^2+ |\nabla u_T|^2}\int_{\R^d}  \tilde \eta_{T,C}.
\end{multline*}
The combination of these three estimates then yields for all $C$ large enough
\begin{multline*}
\Big|\expec{\frac1Tu_T v \mathds1_{\R^d\setminus \calB}+ \nabla u_T \cdot \nabla v
-v( \bar g_\theta \mathds 1_{\R^d\setminus\calB}-\bar g \mathds 1_{\calB})}\Big|\\
\lesssim\,  \frac1{C\sqrt{T}} \expec{\frac1Tu_T^2+ v^2+ |\nabla u_T|^2+|\nabla v|^2+\bar g^2}.
\end{multline*}
The claimed weak formulation \eqref{e.eq-weak-prob} of the equation then follows by letting the parameter $C$ of the cut-off go to infinity.

\subsection{Lemma~\ref{lem:compactness}: Compactness result}
We split the proof into three steps. 
In the first step we derive a reverse Poincar\'e inequality on $B^{R}$ which is at the origin of the compactness result
of the second step that relies on the spectral decomposition of the Laplacian on $B^{2R}$. We then conclude 
in the last step by stationarity of $\nabla u_T$.

\medskip

\step1 Reverse Poincar\'e inequality: For all $R \gg 1$, 
\begin{equation}\label{e.reverse-Poincare}
\int_{B^{R}}   |\nabla u_T|^2 \,\le\,\frac1{R^{2}} \inf_{t\in \R} \int_{B^{2R}} (u_T- t)^2 + C R^2 \int_{B^{2R}} (\frac{1}{T^2}u_T^2+1).
\end{equation}
Let $\eta_R$ be a smooth cut-off for $B^R$ in $B^{2R}$ chosen in a such a way that
\begin{itemize}
\item for any inclusion $B_i$
that intersects $B^{2R}$ the cut-off $\eta_R$ is constant on $\bar B_i$ (where $\bar B_i$ is a ball
with the same center as $B_i$ and twice its radius; $\bar \calB=\cup_i \bar B_i$),
\item for all $\bar B\in \bar \calB$ such that $\bar B\cap \partial B^{2R}\neq \emptyset$ we have $\eta_R|_{\bar B}\equiv 0$,
\item  $|\nabla \eta_R|\lesssim R^{-1}$.
\end{itemize}
Such a cut-off exists for $R\gg 1$: it suffices to consider the standard cut-off $\tilde \eta_R$ for $B^R$ in $B^{2R}$, and modify it locally on each $\bar B_i$ as follows.
In $B_i$ we set $\eta_R \equiv \fint_{\bar B_i} \tilde \eta_R$, on $\partial \bar B_i$ we set $\eta_R \equiv \tilde \eta_R$, and we extend $\eta_R$ radially and linearly
between $\partial B_i$ and $\partial \bar B_i$. So defined, we still have $|\nabla \eta_R|\lesssim R^{-1}$ in the enlarged inclusions $\bar B_i$.

\medskip

We then proceed to the Caccioppoli argument, and consider the test-function $\eta_R^2 (u_T-t)$ for $t\in \R$. Let $I=\{i\in \N\,|\, \bar B_i \subset B^{2R}\}$ and for all $i \in I$ set $u_{T,i}:=u_T |_{B_i}$ and $\eta_{R,i} = \eta_R|_{B_i}$.
This yields
$$
\int_{\R^d \setminus \calB} \eta_R^2 (u_T-t)  \bar g_\theta =\frac1T \int_{\R^d \setminus \calB} \eta_R^2 (u_T-t)u_T+ \int \nabla (\eta_R^2 (u_T-t)) \cdot \nabla u_T + \sum_{i\in I} \bar g \eta_{R,i}^2 (u_{T,i}-t),
$$
that can be rewritten by the properties of $\eta_R$ as
\begin{eqnarray*}
\lefteqn{\frac1T \int_{\R^d \setminus \calB} \eta_R^2 (u_T-t)^2+\int |\nabla (\eta_R (u_T-t))|^2} 
\\
&=& \int (u_T-t)^2 |\nabla \eta_R|^2 - \sum_{i\in I}\bar g \eta_{R,i}^2 (u_{T,i}-t)+\int_{\R^d\setminus \calB} \eta_R^2(u_T-t)\bar g_\theta
-\frac1T t \int_{\R^d \setminus \calB}\eta_R^2 (u_T-t)
\\
&\lesssim & \frac1{R^2} \int_{B^{2R}} |u_T-t|^2+\int_{B^{2R}} |u_T-t|+\frac1T |t| \int_{B^{2R}} |u_T-t|.
\end{eqnarray*}
By Young's inequality on the last two RHS terms in form of $|ab|\le \frac1{2CR^2} a^2+\frac{CR^2}2 b^2$,
and the choice $t=\fint_{B^{2R}} u_T$ that entails the bound $t^2 \,\le \, \fint_{B^{2R}} u_T^2$,  this implies the desired estimate \eqref{e.reverse-Poincare} for a suitable $C\sim 1$.

\medskip

\step2 Proof of \eqref{e.compact0}.

This proof is similar to \cite{Gloria-Otto-10b}; it is reproduced for completeness.
In view of \eqref{e.reverse-Poincare}, it is enough to prove that for any function $v\in H^1(B^{2R})$,
\begin{equation}\label{e.compact0bis}
\frac1{R^{2}} \inf_{t\in \R} \int_{B^{2R}} (v- t)^2 \le C \sum_{n=1}^N \Big(\int_{B^{2R}}  \nabla v\cdot R^{-\frac d2} g_n(\tfrac \cdot R)\Big)^2 + \delta \int_{B^{2R}} |\nabla v|^2.
\end{equation}
By rescaling length according to $x=R\hat x$,
we may assume that $2R=1$. Let $\{(\lambda_n,u_n)\}_{n=0,1,\cdots}$ denote 
a complete set of increasing eigenvalues and $L^2$-orthonormal eigenfunctions of $-\triangle$ on $B^1$ endowed
with homogeneous Neumann boundary conditions, that is
\begin{equation}\label{eq:co-u.18}
\int_{B^1}\nabla v\cdot\nabla u_ndx=\lambda_n\int_{B^1}v u_ndx\quad
\mbox{for all functions }\;v\in H^1(B^1).
\end{equation}
In particular, we have $\int_{B^1}|\nabla u_n|^2dx=\lambda_n\int_{B^1}u_n^2dx=\lambda_n$.
We also note that $\lambda_1>0$. Hence for all $n\ge 1$ 
\begin{equation}\label{eq:co-u.19}
F_nu=\int_{B^1}\nabla u\cdot\frac{\nabla u_n}{\sqrt{\lambda_n}}dx\quad\mbox{for all functions}\;u\in H^1(B^1)
\end{equation}
defines a linear functional $F_n$ on vector fields that has the boundedness property $(F_nu)^2\le\int_{B^1}|\nabla u|^2dx$.
By completeness of the orthonormal system $\{u_n\}_{n=0,1,\cdots}$, Plancherel
and with $\bar u=\fint_{B^1} u$, we have
\begin{eqnarray*}
\lefteqn{\int_{B^1}(u-\bar u)^2dx=\sum_{n=1}^\infty\left(\int_{B^1}u u_ndx\right)^2}\\
&\stackrel{(\ref{eq:co-u.18})}{=}&
\sum_{n=1}^\infty\frac{1}{\lambda_n}
\left(\int_{B^1}\nabla u\cdot \frac{\nabla u_n}{\sqrt{\lambda_n}}dx\right)^2\\
&\le&\frac{1}{\lambda_1}\sum_{n=1}^{N-1}\left(\int_{B^1}\nabla u\cdot \frac{\nabla u_n}{\sqrt{\lambda_n}}dx\right)^2
+\frac{1}{\lambda_{N}}\sum_{n=N}^\infty\left(\int_{B^1}\nabla u\cdot \frac{\nabla u_n}{\sqrt{\lambda_n}}dx\right)^2.
\end{eqnarray*}
We note that $\eqref{eq:co-u.18}$ yields that also $\{\frac{\nabla u_n}{\sqrt{\lambda_n}}\}_{n=1,\cdots}$
is orthonormal, so that the above together with definition \eqref{eq:co-u.19} yields
\begin{eqnarray*}
\int_{B^1}(u-\bar u)^2dx
&\le&\frac{1}{\lambda_1}\sum_{n=1}^{N-1}(F_n\nabla u)^2
+\frac{1}{\lambda_{N}}\int_{B^1}|\nabla u|^2dx.
\end{eqnarray*}
Because of $\lim_{N\uparrow\infty}\lambda_{N}=\infty$, this implies \eqref{e.compact0bis} in its
$(2R=1)$-version, and therefore  \eqref{e.compact0}.

\medskip

\step3 Proof of \eqref{e.compact}.

With the choice $\delta = 2^{-d-2}$, the expectation of \eqref{e.compact0} to the power $p\ge 1$ yields 
\begin{multline*}
\expec{\Big(\int_{B^R} |\nabla u_T|^2\Big)^p}^\frac1p \,\le \, C \sum_{n=1}^N \expec{\Big(\int_{B^{2R}}  \nabla u_T\cdot g_n\Big)^{2p}}^\frac1p
\\
+ 2^{-d-2} \expec{ \Big(\int_{B^{2R}} |\nabla u_T|^2\Big)^p}^\frac1p+ CR^2\Big(1+\frac1{T^2} \expec{\Big(\int_{B^{2R}} u_T^2\Big)^p}^\frac1p\Big).
\end{multline*}
We then use Poincar\'e's inequality on the last RHS term to the effect of 
$$
\frac{R^2}{T^2} \expec{\Big(\int_{B^{2R}} u_T^2\Big)^p}^\frac1p \,\le \, C\frac{R^4}{T^2} \expec{\Big(\int_{B^{2R}} |\nabla u_T|^2\Big)^p}^\frac1p
+ C\frac{R^2}{T} \expec{\Big(\int_{B^{2R}} \frac1{\sqrt T}u_T\Big)^{2p}}^\frac1p.
$$
Provided $\sqrt T\gg R$, the desired estimate \eqref{e.compact} then follows from absorbing the terms $\expec{\Big(\int_{B^{2R}} |\nabla u_T|^2\Big)^p}^\frac1p$ 
into the LHS by stationarity of $\nabla u_T$.

\subsection{Lemma~\ref{lem:diff-u}: Equation for differences}

In the following calculations, we assume that $w_T \in H^1(\R^d)$ -- which can be checked a posteriori.
In particular, testing \eqref{e.diff-uT} with $w_T$ itself yields
\begin{eqnarray*}
\int_{\R^d} \frac1T w_T^2+|\nabla w_T|^2&=& \bar g_\theta \int_{\R^d} w_T(\mathds 1_{\R^d\setminus \calB'}-\mathds 1_{\R^d\setminus \calB})+\frac1T \int_{\R^d} w_T(u'_T \mathds1_{\calB'}-u_T \mathds1_{\calB})
\\
&&-\sum_i \int_{\partial B_i'} w_T \nabla u'_T \cdot n'_i+\int_{\partial B_i} w_T \nabla u_T \cdot n_i .
\end{eqnarray*}
By the definition of  $\calB\triangle \calB'$   we have
\begin{eqnarray*}
\Big| \int_{\R^d} w_T(\mathds 1_{\R^d\setminus \calB'}-\mathds 1_{\R^d\setminus \calB}) \Big| &\le& \int_{\calB\triangle \calB'} |w_T|,
\\
\Big|\int_{\R^d} w_T(u'_T \mathds1_{\calB'}-u_T \mathds1_{\calB})\Big|&\le& \int_{\calB'} w_T^2 +\int_{\calB\triangle \calB'} |w_T||u_T|,
\end{eqnarray*}
and since $w_T$ is constant on balls $B \in \calB \cap \calB'$ and $\int_{\partial B_i} \nabla u_T\cdot n_i=1$ (likewise for $u_T'$)
\begin{multline*}
\Big|\sum_i \int_{\partial B_i'} w_T \nabla u'_T \cdot n'_i- \int_{\partial B_i} w_T \nabla u_T \cdot n_i\Big|
\\
\le\, \sum_{B_i \in \calB \setminus \calB'} \Big|\int_{\partial B_i} w_T \nabla u_T \cdot n_i\Big|+ \sum_{B_i' \in \calB' \setminus \calB} \Big|\int_{\partial B_i'} w_T \nabla u'_T \cdot n_i'\Big|.
\end{multline*}
The combination of these four estimates yields
\begin{multline}\label{e.pr-diff-uT}
\int_{\R^d \setminus \cup_{B'\in \calB'}B'} \frac1T w_T^2+\int_{\R^d} |\nabla w_T|^2\,\lesssim\, \int_{\calB\triangle \calB'} |w_T|
+\frac1T \int_{\calB\triangle \calB'} |w_T||u_T|
\\
+\sum_{B_i \in \calB \setminus \calB'} \Big|\int_{\partial B_i} w_T \nabla u_T \cdot n_i\Big|+ \sum_{B_i' \in \calB' \setminus \calB} \Big|\int_{\partial B_i'} w_T \nabla u'_T \cdot n_i'\Big|.
\end{multline}
We first treat the first two RHS terms of \eqref{e.pr-diff-uT}.
By H\"older's inequality, the Poincar\'e-Sobolev inequality for $d>2$, and Young's inequality with constant $C\gg1$ (to be fixed later)
\begin{eqnarray*}
 \int_{\calB\triangle \calB'} |w_T|&\le &  |\calB\triangle \calB'|^\frac{d+2}{2d}\Big(\int_{\R^d} |w_T|^\frac{2d}{d-2}\Big)^\frac{d-2}{2d}
\\
&\lesssim&   |\calB\triangle \calB'|^\frac{d+2}{2d}\Big(\int_{\R^d} |\nabla w_T|^2\Big)^\frac12
\\
&\lesssim& \frac1C \int_{\R^d} |\nabla w_T|^2+C  |\calB\triangle \calB'|^\frac{d+2}{d}.
\end{eqnarray*}
Likewise, we obtain for the second RHS term
\begin{eqnarray*}
\frac1T \int_{\calB\triangle \calB'} |w_T||u_T|&\lesssim& \frac1C \int_{\R^d} |\nabla w_T|^2+C \frac1{T^2} \Big(\int_{\calB\triangle \calB'} |u_T|^\frac{2d}{d+2}\Big)^\frac{d+2}{d}
\\
&\lesssim&  \frac1C \int_{\R^d} |\nabla w_T|^2+C  |\calB\triangle \calB'|^{\frac 2d} \frac1{T^2} \int_{\calB\triangle \calB'} u_T^2.
\end{eqnarray*}
It remains to estimate the last RHS sums of \eqref{e.pr-diff-uT}.
Let $B_i \in \calB\setminus \calB'$ and set $\bar w_{T,i}=\fint_{B_i} w_T$ and $\tilde g=\frac{1}{|\partial B|}$.
Since $\int_{\partial B_i} \nabla u_T \cdot n_i=1$, we have
\begin{eqnarray*}
\int_{\partial B_i} w_T \nabla u_T \cdot n_i&=&\int_{\partial B_i} w_T \tilde g+\int_{\partial B_i} w_T (\nabla u_T \cdot n_i-\tilde g)
\\
&=&\int_{\partial B_i} w_T \tilde g+\int_{\partial B_i} (w_T-\bar w_{T,i}) (\nabla u_T \cdot n_i-\tilde g).
\end{eqnarray*}
By elliptic regularity up to the boundary for both $u_T$ and $u'_T$ taken separately, $\sup_{B_i} |w_T - \bar w_{T,i}| \,\lesssim \, \Big(\int_{\bar B_i} |\nabla u_T|^2+ |\nabla u'_T|^2\Big)^\frac12\,\lesssim\, \Big(\int_{\bar B_i} |\nabla w_T|^2+ |\nabla u_T|^2\Big)^\frac12$, so that
\begin{eqnarray*}
\lefteqn{\Big| \int_{\partial B_i} w_T \nabla u_T \cdot n_i\Big| }
\\
&\lesssim& \int_{B_i} |w_T|+\Big(\int_{\bar B_i} |\nabla w_T|^2+ |\nabla u_T|^2\Big)^\frac12\Big(1+\int_{\partial B_i} |\nabla u_T \cdot n_i-\tilde g|\Big)
\\
&\lesssim &  \int_{B_i} |w_T| +\frac1C \int_{\bar B_i}( |\nabla w_T|^2+|\nabla u_T|^2)+C\Big(1+\int_{\partial B_i} |\nabla u_T \cdot n_i-\tilde g|\Big)^2.
\end{eqnarray*}
Likewise, for  $B_i' \in \calB'\setminus \calB$, we have
\begin{eqnarray*}
\Big| \int_{\partial B_i'} w_T \nabla u_T' \cdot n_i\Big| &\lesssim & \int_{B_i'} |w_T| +\frac1C \int_{\bar B_i'}( |\nabla w_T|^2+|\nabla u_T|^2)+C\Big(1+\int_{\partial B_i'} |\nabla u_T' \cdot n_i'-\tilde g|\Big)^2.
\end{eqnarray*}
Summing these two estimates over $i$ then yields
\begin{multline*}
{\sum_{B_i \in \calB\setminus \calB'} \Big|\int_{\partial B_i} w_T \nabla u_T \cdot n_i\Big|+ \sum_{B_i' \in \calB'\setminus \calB} \Big|\int_{\partial B_i'} w_T \nabla u'_T \cdot n_i'\Big|}
\\
\,\lesssim\, |\bar g|\int_{\calB\triangle \calB'}  |w_T|+\frac1C \int_{\R^d} |\nabla w_T|^2+C|\calB\triangle \calB'|
\\
+\sum_{B_i \in \calB\setminus \calB'} \Big(C\Big(\int_{\partial B_i} |\nabla u_T \cdot n_i-\tilde g|\Big)^2+\frac1C \int_{\bar B_i}|\nabla u_T|^2\Big)
\\
+\sum_{B_i' \in \calB'\setminus \calB} \Big(C\Big(\int_{\partial B_i'} |\nabla u_T' \cdot n_i'-\tilde g|\Big)^2+\frac1C \int_{\bar B_i'}|\nabla u_T|^2\Big)
.
\end{multline*}
As for the first RHS term of \eqref{e.pr-diff-uT}, we have
$$
|\bar g| \int_{\calB\triangle \calB'} |w_T|\,\lesssim\, \frac1C \int_{\R^d} |\nabla w_T|^2+C\bar g^2  |\calB\triangle \calB'|^\frac{d+2}{d}.
$$
The combination of all these estimates for some $C\gg 1$ large enough to absorb terms of the form $ \int_{\R^d} |\nabla w_T|^2$ into the LHS
finally gives
\begin{multline*}
\int_{\R^d\setminus \cup_{B'\in \calB'} B'}\frac1T w_T^2+\int_{\R^d} |\nabla w_T|^2 \,\le\,  C(\bar g^2+\bar g^2_\theta)  |\calB\triangle \calB'|^\frac{d+2}{d}+C  |\calB\triangle \calB'|^{\frac 2d} \frac1{T^2} \int_{\calB\triangle \calB'} u_T^2
\\
+C\int_{\overline D} |\nabla u_T|^2
+\sum_{B_i \in \calB \setminus  \calB'} C\Big(\int_{\partial B_i} |\nabla u_T \cdot n_i-\tilde g|\Big)^2
+\sum_{B_i' \in \calB' \setminus \calB} C\Big(\int_{\partial B_i'} |\nabla u_T' \cdot n_i'-\tilde g|\Big)^2.
\end{multline*}
The desired bound \eqref{e.diff-uT-estim} now follows in combination with the trace estimate $\int_{B'} w_T^2 \lesssim \int_{\overline{B'}\setminus B'} w_T^2+|\nabla w_T|^2$ on each ball $B'\in \calB'$.

\subsection{Corollary~\ref{cor:apriori}: Estimates of differences}
Let $\rho_0\ge \rho_\m$.
By Lemma~\ref{lem:diff-u}, it suffices to control the terms 
$$\sum_{B_i \in \calB\setminus \calB'} C\Big(\int_{\partial B_i} |\nabla u_T \cdot n_i-\tilde g|\Big)^2+\sum_{B_i' \in \calB'\setminus \calB} C\Big(\int_{\partial B_i'} |\nabla u_T' \cdot n_i'-\tilde g|\Big)^2.$$
We apply Lemma~\ref{lem:reg1} to $u_T$ on balls $B_i$ and use in addition a trace estimate on $\partial B_i$ to control $\bar v^2$,  to the effect of
$$
\Big(\int_{\partial B_i} |\nabla u_T \cdot n_i-\tilde g|\Big)^2 \,\le \, \rho_0^{-2\alpha} \Big(\int_{B_i^{\rho_0}}(|\nabla u_T|^2+\frac{\rho_0^{d+4}}{T^2}(|\nabla u_T|^2+u_T^2))+\rho_0^{d} (1+\frac{\rho_0^2}T)(1+\rho_0^2)\Big),
$$
which we write as
$$
\Big(\int_{\partial B_i} |\nabla u_T \cdot n_i-\tilde g|\Big)^2 \,\le \, 2\rho_0^{-2\alpha}  \Big(\rho_0^{d+4}\int_{B_i^{\rho_0}}( |\nabla u_T|^2+\frac1{T^2}  u_T^2)+\rho_0^{d} (1+\rho_0^2 )\Big).
$$
Likewise, for $u_T'$ this yields by the triangle inequality
\begin{multline*}
\Big(\int_{\partial B_i'} |\nabla u_T' \cdot n_i'-\tilde g|\Big)^2
\,
\le \,4 \rho_0^{-2\alpha}  \rho_0^{d+4} \int_{(B'_i)^{\rho_0}} (|\nabla u_T|^2+\frac1{T^2}u_T^2)+  \\
4\rho_0^{-2\alpha} \int_{(B'_i)^{\rho_0}}(|\nabla w_T|^2+\frac{\rho_0^{d+4}}{T^2}(|\nabla w_T|^2+w_T^2))
+2 \rho_0^{-2\alpha} \rho_0^{d} (1+ \rho_0^2).
\end{multline*}
We then choose $\rho_0>0$ and $T\gg \rho_0^{d+4}$ so large that $4\rho_0^{-2\alpha}(1+\frac{\rho_0^{d+4}}{T})C\le \frac14$, and then conclude by absorbing the term $\sum_{B_i' \in \calB'\setminus \calB} 2 \rho_0^{-2\alpha}C \int_{(B'_i)^{\rho_0}}(1+ \frac{\rho_0^{d+4}}{T})(|\nabla w_T|^2+\frac1{T}w_T^2)\le \frac14 \int_{\R^d} (|\nabla w_T|^2+\frac1{T}w_T^2)$
into the LHS of~\eqref{e.diff-uT-estim}.

\section*{Acknowledgements}
The author warmly thanks Mitia Duerinckx, Arianna Giunti,
Richard H\"ofer, Jonas Jansen, Jules Pertinand, and Juan Vel\'azquez for 
inspiring discussions on sedimentation, and on the scalar model studied here,
as well as David G\'erard-Varet and Sylvia Serfaty for pointing out the relation to 
the Coulomb energy of point processes.
Financial support is acknowledged from the European Research Council under the European Community's Seventh Framework Programme (FP7/2014-2019 Grant Agreement QUANTHOM 335410).

\bibliographystyle{plain}
%\bibliography{biblio}

\begin{thebibliography}{10}

\bibitem{sedimentation}
{Sedimentation -- Learning \& Informative blog for Engineering Students}, 2012.
\newblock \href{http://engrsl.blogspot.fr/2012/04/sedimentation.html}{Webpage}.

\bibitem{AKM2}
S.~N. {Armstrong}, T.~{Kuusi}, and J.-C. {Mourrat}.
\newblock The additive structure of elliptic homogenization.
\newblock {\em Invent. Math.}, 208:999--1154, 2017.

\bibitem{AKMbook}
S.~N. {Armstrong}, T.~{Kuusi}, and J.-C. {Mourrat}.
\newblock {\em Quantitative Stochastic Homogenization and Large-Scale
  Regularity}.
\newblock Grundlehren der mathematischen Wissenschaften. Springer International
  Publishing, 2019.

\bibitem{AS}
S.~N. Armstrong and C.~K. Smart.
\newblock Quantitative stochastic homogenization of convex integral
  functionals.
\newblock {\em Ann. Sci. \'Ec. Norm. Sup\'er. (4)}, 49(2):423--481, 2016.

\bibitem{Batchelor}
G.~K. Batchelor.
\newblock Sedimentation in a dilute dispersion of spheres.
\newblock {\em J. Fluids Mech.}, 52(2):245--268, 1972.

\bibitem{Blanc-Lewin-15}
X.~Blanc and M.~Lewin.
\newblock The crystallization conjecture: a review.
\newblock {\em EMS Surv. Math. Sci.}, 2(2):225--306, 2015.

\bibitem{pete}
P.~Bower.
\newblock {Pete's Lab: Sediment Stirring, Columbia University}, 2011.
\newblock \href{https://www.youtube.com/watch?v=JI-xah5O6vk}{Movie}.

\bibitem{MR0014916}
J.~M. Burgers.
\newblock On the influence of the concentration of a suspension upon the
  sedimentation velocity (in particular for a suspension of spherical
  particles).
\newblock {\em Nederl. Akad. Wetensch., Proc.}, 44:1045--1051, 1177--1184; 45,
  9--16, 126--128 (1942), 1941.

\bibitem{MR815929}
R.~E. Caflisch.
\newblock Sedimentation of a random dilute suspension.
\newblock In {\em Macroscopic modelling of turbulent flows ({N}ice, 1984)},
  volume 230 of {\em Lecture Notes in Phys.}, pages 14--23. Springer, Berlin,
  1985.

\bibitem{Caflisch-Luke}
R.E. Caflisch and J.H.C. Luke.
\newblock Variance in the sedimentation speed of a suspension.
\newblock {\em The Physics of fluids}, 28(3):759--760, 1985.

\bibitem{Courant-Hilbert-II}
R.~Courant and D.~Hilbert.
\newblock {\em Methods of mathematical physics. {V}ol. {II}}.
\newblock Wiley Classics Library. John Wiley \& Sons, Inc., New York, 1989.
\newblock Partial differential equations, Reprint of the 1962 original, A
  Wiley-Interscience Publication.

\bibitem{DG2}
M.~Duerinckx and A.~Gloria.
\newblock {Multiscale functional inequalities in probability: Concentration
  properties}.
\newblock {\em ALEA, Lat. Am. J. Probab. Math. Stat.}, 17:133--157, 2020.

\bibitem{DG1}
M.~Duerinckx and A.~Gloria.
\newblock {Multiscale functional inequalities in probability: Constructive
  approach}.
\newblock {\em Annales Henri Lebesgue}, 2020.
\newblock In press.

\bibitem{GVH-19}
D.~{G\'erard-Varet} and M.~{Hillairet}.
\newblock {Analysis of the viscosity of dilute suspensions beyond Einstein's
  formula}.
\newblock {\em arXiv e-prints}, page arXiv:1905.08208, May 2019.

\bibitem{GNO-quant}
A.~Gloria, S.~Neukamm, and F.~Otto.
\newblock Quantitative estimates in stochastic homogenization for correlated
  coefficient fields.
\newblock Preprint, arXiv:1910.05530.

\bibitem{GNO1}
A.~Gloria, S.~Neukamm, and F.~Otto.
\newblock Quantification of ergodicity in stochastic homogenization: optimal
  bounds via spectral gap on {G}lauber dynamics.
\newblock {\em Invent. Math.}, 199(2):455--515, 2015.

\bibitem{GNO-reg}
A.~Gloria, S.~Neukamm, and F.~Otto.
\newblock A regularity theory for random elliptic operators.
\newblock {\em Milan Journ. Mathematics}, 2020.
\newblock In press.

\bibitem{GO4}
A.~Gloria and F.~Otto.
\newblock The corrector in stochastic homogenization: optimal rates, stochastic
  integrability, and fluctuations.
\newblock Preprint, arXiv:1510.08290.

\bibitem{GO1}
A.~Gloria and F.~Otto.
\newblock An optimal variance estimate in stochastic homogenization of discrete
  elliptic equations.
\newblock {\em Ann. Probab.}, 39(3):779--856, 2011.

\bibitem{Gloria-Otto-10b}
A.~Gloria and F.~Otto.
\newblock Quantitative results on the corrector equation in stochastic
  homogenization.
\newblock {\em Journ. Europ. Math. Soc. (JEMS)}, 19:3489--3548, 2017.

\bibitem{MR2768013}
E.~Guazzelli and J.~Hinch.
\newblock Fluctuations and instability in sedimentation.
\newblock In {\em Annual review of fluid mechanics. {V}olume 43, 2011},
  volume~43 of {\em Annu. Rev. Fluid Mech.}, pages 97--116. Annual Reviews,
  Palo Alto, CA, 2011.

\bibitem{HM-12}
B.~M. Haines and A.~L. Mazzucato.
\newblock A proof of {E}instein's effective viscosity for a dilute suspension
  of spheres.
\newblock {\em SIAM J. Math. Anal.}, 44(3):2120--2145, 2012.

\bibitem{Hofer-18}
R.~M. H\"{o}fer.
\newblock Sedimentation of inertialess particles in {S}tokes flows.
\newblock {\em Comm. Math. Phys.}, 360(1):55--101, 2018.

\bibitem{JO-04}
P.-E. Jabin and F.~Otto.
\newblock Identification of the dilute regime in particle sedimentation.
\newblock {\em Comm. Math. Phys.}, 250(2):415--432, 2004.

\bibitem{Leble}
T.~Lebl\'{e}.
\newblock Logarithmic, {C}oulomb and {R}iesz energy of point processes.
\newblock {\em J. Stat. Phys.}, 162(4):887--923, 2016.

\bibitem{Penrose-01}
M.~D. Penrose.
\newblock Random parking, sequential adsorption, and the jamming limit.
\newblock {\em Comm. Math. Phys.}, 218(1):153--176, 2001.

\bibitem{Petrache-Serfaty-19}
M.~{Petrache} and S.~{Serfaty}.
\newblock {Crystallization for Coulomb and Riesz Interactions as a Consequence
  of the Cohn-Kumar Conjecture}.
\newblock {\em arXiv e-prints}, page arXiv:1908.09714, Aug 2019.

\bibitem{MR1021642}
J.~Rubinstein and J.~B. Keller.
\newblock Sedimentation of a dilute suspension.
\newblock {\em Phys. Fluids A}, 1(4):637--643, 1989.

\bibitem{Sandier-Serfaty-15}
E.~Sandier and S.~Serfaty.
\newblock 2{D} {C}oulomb gases and the renormalized energy.
\newblock {\em Ann. Probab.}, 43(4):2026--2083, 2015.

\bibitem{Serfaty-14}
S.~{Serfaty}.
\newblock {\em {Coulomb Gases and Ginzburg-Landau Vortices}}, volume~21.
\newblock EMS Zurich Lectures in Advanced Mathematics, 2015.

\bibitem{Theil-06}
F.~Theil.
\newblock A proof of crystallization in two dimensions.
\newblock {\em Comm. Math. Phys.}, 262(1):209--236, 2006.

\end{thebibliography}

\def\cprime{$'$} \def\cprime{$'$} \def\cprime{$'$}

\end{document}